 \documentclass[final,1p,times,review]{elsarticle}

\usepackage{amssymb}
\usepackage{mathbbold,bbm,bm}
\usepackage{mathrsfs,amsmath,amssymb,booktabs,array,xcolor,url,color,soul,multirow,multicol,balance}
\usepackage[justification=centering]{caption}
\usepackage{graphicx}
\usepackage{subfigure}
\usepackage[all]{xy}
\usepackage{chemarrow}

\journal{System \& Control Letters}

\newtheorem{theorem}{Theorem}[section]
\newtheorem{definition}[theorem]{Definition} 
\newtheorem{remark}[theorem]{Remark} 
\newtheorem{lemma}[theorem]{Lemma}
\newtheorem{corollary}[theorem]{Corollary}
\newtheorem{proposition}[theorem]{Proposition}
\begin{document}

\begin{frontmatter}

\title{Revisiting Persistence of Chemical Reaction Networks through Lyapunov Function Partial Differential Equations\tnoteref{label1}}
\tnotetext[label1]{This work is supported by the National Nature Science Foundation of China under [Grant No. 11671418, 61611130124].}

 \author[label1]{Xiaoyu~Zhang}
 \ead{11735035@zju.edu.cn}
\author[label2]{ Zhou~Fang}
\ead{zhou\_fang@zju.edu.cn}
 \author[label3]{Chuanhou~Gao \corref{cor1}}
 \ead{gaochou@zju.edu.cn}
 \cortext[cor1]{Corresponding author}
\address[label1,label2,label3]{School of Mathematical Sciences, Zhejiang University, Hangzhou 310027, China}



\begin{abstract}
In this paper, we propose a novel technique, referred to as the Lyapunov Function PDEs technique, to diagnose persistence of chemical reaction networks with mass-action kinetics. The technique allows that every network is attached to a Lyapuonv function PDE and a boundary condition whose solutions are expected to be Lyapunov functions for the network. By means of solution of the PDEs, either in the forms of itself or its time derivative, some checkable criteria are proposed for persistence of network systems. These criteria show high validity in proving that neither non-semilocking boundary points nor semilocking boundary non-equilibrium points (with additional conditions included) in mass-action systems are $\omega$-$limit$ points. Further, we prove that a class of networks, called $\mathcal{W}_{I}$-endotactic networks that are regardless of values of the system parameters, a set larger than endotactic networks but smaller than $w_{I}$-endotactic networks, also give rise to persistent systems if the networks are $1$-dimensional. Although part of our results are covered by existing ones, part of them are still new. The proposed Lyapunov Function PDEs technique allows us to check persistence of mass action systems in an alternative corner, and is exhibiting large potential. 
\end{abstract}

\begin{keyword}
chemical reaction networks\sep mass-action kinetics\sep persistence\sep Lyapunov Function PDEs\sep $\omega$-$limit$ point. 
\end{keyword}

\end{frontmatter}


\section{Introduction}
Persistence is a concept from population processes that can be interpreted as non-extinction. A system of species interactions is said to be persistent if initial present species are not allowed to approach or asymptotically approach extinction. There are many practical examples that can be modeled by population processes \cite{Anderson2010}, like animal populations, the spread of infectious diseases, and the evolution of biochemical reaction systems. The species in these examples might correspond to types of animals, infected and noninfected individuals, and chemical reactants and resultants, respectively. Persistence is thus a property of central practical interest. Also, it is of major theoretical interest since from a mathematical perspective persistence means that the $\omega$-$limit$ set of any trajectory starting in the interior of the positive orthant does not intersect the boundary of the positive orthant \cite{Angeli2007}. In this work, we focus on the property of persistence for mass-action systems (MASs), i.e., chemical reaction networks (CRNs) equipped with mass-action kinetics, and revisit some known results about persistence except something new using a novel technique; that is, Lyapunov Function partial differential equations (PDEs) \cite{Fang}. 

Mass-action CRNs often give rise to a family of nonequilibrium ordinary differential equation (ODE) systems that lie at the limits of what we can treat with the current tools of dynamical systems. Although the past $40$ years have witnessed the fast development of chemical reaction network theory (CRNT) since the pioneering work of Horn and Jackson \cite{Horn1972,Horn1974}, and Feinberg (see, e.g., \cite{Feinberg1987,Feinberg1995}), it still poses a great challenge to determine qualitative property of solutions of mass-action systems. Several basic questions about reaction networks remain open, chief among them being the global attractor conjecture (GAC).

\textit{Global Attractor Conjecture. A complex balanced equilibrium contained in the interior of a positive stoichiometric compatibility class is is a global attrator of the interior of that compatibility class.}

This conjecture was initiated by Horn \cite{Horn1974} in $1974$, and named GAC by Craciun et. al. \cite{Craciun2009}. A more general conjecture is persistence conjecture presented by Feinberg \cite{Feinberg1987} that says any weakly reversible mass-action system to be persistent. Here, a weakly reversible mass-action system is that each of the connected components of its directed reaction graph is strongly connected (see \textbf{Definition \ref{def:2reversible}}). Since the weakly reversible structure is a necessary network requirement for complex balanced mass-action systems, that persistence conjecture is true might lead to GAC to hold. A further result is that persistence coupled with complex balancing will guarantee global convergence to equilibria in complex balanced mass-action systems. Partly due to this reason, and partly due to universal presence of CRNs in systems biology and bioengineering, there has been much activity in recent years aimed at the resolution of these two conjectures (see, e.g., \cite{Anderson2010,Anderson2011,Craciun2013}). The studies related to persistence even go beyond the structure of weakly reversible mass-action systems \cite{Angeli2007,Angeli2011,Pantea2012,G2014}. Some new concepts, like semilocking set, face, endotatic network (see \textbf{Definition \ref{def:semilocking}, \ref{def:face}, \ref{def:1endotatic}, \ref{def:Endotatic}}), etc., were accordingly proposed for persistence checking. 

The meaning of persistence in CRNT is analogous to that in population processes. Namely, none of chemical species in a MAS can be completely "used up" in the course of the reaction if every species has non-zero initial concentration. Mathematically, it measures if any trajectory in a MAS will tend to the face of $\mathbb{R}^{n}_{\geq 0}$. Based on the notions of semilocking set, face and endotactic network, many efforts have been made to characterize persistence of CRNs. Angeli et al. \cite{Angeli2007} proved that if a boundary point (an element in $L_W$) is the $\omega$-$limit$ point, then the species set $W$ in which every species has zero concentration must be a semilocking set. They also reported that the global conservative CRNs with each semilocking set having a  conservation relation between some species in it are persistent. In addition, they obtained that the global conservative CRNs with some semilocking sets without local conservation relation in it are persistent if the following two conditions are true: (a) the semilocking sets are dynamically non-emptiable; (b) there are no nested distinct locking sets without local conservation relation in the network. Johnston and Siegel \cite{Johnson2011} extended this result by putting forward that bounded networks with all weak dynamical non-emptiable semilocking sets have persistence. Anderson and his coauthors \cite{Anderson2008,Anderson2010,Anderson2011} made a systematic study on GAC. They said that (a) deficiency zero weakly reversible MASs have persistence if every boundary equilibrium in $L_{W}$ ($W$ is a semilocking set in the system) is isolated \cite{Anderson2008}; (b) complex balanced MASs are persistent if the equilibria are in the relative interior of $F_{W}$ (a codimension-one face compared to the stoichiometric compatiblity class, see details in \textbf{Definition \ref{def:face}} and the following) \cite{Anderson2010}; (c) all complex-balanced networks with single linkage class are persistent\cite{Anderson2011}. A recently proposed proof of this conjecture by Craciun \cite{Craciun2015} is currently under verification. Craciun and his coauthors \cite{Pantea2012,Craciun2013,G2014} also studied persistence of CRNs with endotactic structure. Pantea \cite{Pantea2012} proved any $2d$ endotactic MAS with bounded trajectories is persistent. Craciun et al. \cite{Craciun2013} stated that any endotactic MAS with two species is persistent, and Gopalkrishnan et al. \cite{G2014} further asserted that any strongly endotactic MAS is persistent.

One of the most powerful tools for characterizing dynamical behaviors about reaction network systems is to obtain a Lyapunov function. When the networks are complex balanced, there is a canonical choice of Lyapunov function, i.e., the well-known pseudo-Helmholtz free energy function \cite{Horn1972}. However, in non-complex balancing case, there is no general method to obtain a Lyapunov function, and no understanding of the space of possible Lyapunov functions for a given reaction network. Although Anderson et al. \cite{Anderson2015} showed that the pseudo-Helmholtz free energy function can be derived from an appropriate limit of the stationary distributions for the stochastic mass-action systems, this is not a practical construction method since for almost all reaction network systems, the stationary distributions are prohibitively hard to obtain. Our early work \cite{Fang} followed up this issue, and by connecting the microscopic and the macroscopic level, thermodynamics and potential theory, proposed a PDE and a boundary condition for a mass-action system, referred to as the Lyapuonv Function PDEs. When the PDEs are solved, the solutions are conjectured as Lyapunov functions for the network attached. This conjecture has been proved true when the network is complex balanced, of $1$-dimensional stoichiometric subspace, and a combination of the other cases, respectively. Some special networks with dimension of the stoichiometric subspace more than $2$ also support the conjecture. In this paper, we continue to utilize the function of the Lyapunov Function PDEs to analyze persistence of mass-action systems. The solution of the the PDEs, either in the form of itself or its time derivative, is used to yield checkable criteria for persistence of the attached network systems. These criteria work well for checking non-semilocking boundary points, semilocking boundary non-equilibrium points (with additional conditions included) in mass-actions systems to be not $\omega$-$limit$ points. They are further applied to $\mathcal{W}_{I}$-endotactic networks, a set larger than endotactic networks but smaller than $w_{I}$-endotactic networks, which give rise to persistent systems if the networks are $1$-dimensional. Note that $\mathcal{W}_{I}$-endotactic networks are a large class of networks independently of the values of the reaction rate parameters. Therefore, the current proposed checkable criteria have somewhat commonality despite only $1$-dimensional $\mathcal{W}_{I}$-endotactic network to have been tackled.  

The remainder of this paper is organized as follows. \textbf{Section \ref{sec:2}} introduces some preliminaries about CRN. This is followed by reviewing the definition, some known checkable criteria for persistence, and the Lyapunov Function PDEs method in \textbf{Section \ref{sec:3}}. \textbf{Section \ref{sec:4}} presents the checkable criteria for persistence of CRNs under the framework of the Lyapunov Function PDEs technique, and some concrete examples are also illustrated correspondingly. In \textbf{Section \ref{sec:5}}, a kind of new networks called $\mathcal{W}_{I}$-endotactic networks are defined, and the Lyapunov Function PDEs technique exhibits high validity in diagnosing persistence of $1d~\mathcal{W}_{I}$-endotactic networks. Finally, \textbf{Section \ref{sec:6}} concludes the paper. 
~\\~\\
\noindent{\textbf{Mathematical Notation:}}\\
\rule[1ex]{\columnwidth}{0.8pt}
\begin{description}
\item[\hspace{+0.8em}{$\mathbb{R}^n, \mathbb{R}^n_{\geq 0},\mathbb{R}^n_{>0}$}]: $n$-dimensional real space, non-negative and positive real space, respectively.
\item[\hspace{+0.8em}{$x^{v_{\cdot i}}$}]: $x^{v_{\cdot i}}\triangleq\prod_{j=1}^{d}x_{j}^{v_{ji}}$, where $x,v_{\cdot i}\in\mathbb{R}^{n}$.
	\item[\hspace{+0.8em}{$\mathrm{Ln}(x)$}]: $\mathrm{Ln}(x)\triangleq\left(\ln{x_{1}}, \cdots, \ln{{x}_{n}} \right)^{\top}$, where $x\in\mathbb{R}^{n}_{>{0}}$.
	\item[\hspace{+0.8em}{$\mathscr{C}^{i}(\cdot~;\ast)$}]: The set of $i$th continuous differentiable functions from "$\cdot$" to "$\ast$".
	\item[\hspace{+0.8em}{supp$~x$}]: Support set, defined by $\text{supp}~ x=\{i|x_{i}\neq 0\},~\forall x\in \mathbb{R}^{n}_{\geq 0}$.
	\item[\hspace{+0.8em}{supp$^c~x$}]: The complementary set of support set, defined by $\text{supp}^c~ x=\{i|x_{i}= 0\},~\forall x\in \mathbb{R}^{n}_{\geq 0}$.
\item[\hspace{+0.8em}{$\mathbbold{0}_{n}$}]: An $n$-dimensional vector withe each element to be zero.
\item[\hspace{+0.8em}{$0^{0}$}]: The result is defined by $1$.
\end{description}
\rule[1ex]{\columnwidth}{0.8pt}

\section{Preliminaries on CRN}
\label{sec:2}
In this section, some fundamental concepts about CRN \cite{Horn1972, Feinberg1995} are reviewed.
\begin{definition}[Chemical Reaction Network]\label{def:1}
A CRN is composed of three finite sets:
\begin{itemize}
  \item[\rm{(1)}] a set of species $\mathcal{S}=\bigcup_{j=1}^{n}\{S_{j}\}$,
\item[\rm{(2)}] a set of complexes $\mathcal{C}=\bigcup_{i=1}^r\{v_{\cdot i},v'_{\cdot i}\}$ satisfying $\mathrm{Card}(\mathcal{C})=c$ and $$\bigcup_{i=1,\cdots,r}\{\mathrm{supp} ~{v_{\cdot i}}\cup \mathrm{supp}~{v'_{\cdot i}}\}=\mathcal{S},$$
where $v_{\cdot i},v'_{\cdot i}\in \mathbb{Z}^{\mathcal{S}}_{\geq 0}$ with the $j$th entry of $v_{\cdot i}$ to express the stoichiometric coefficient of species $S_{j}$ in $v_{\cdot i}$,
\item[\rm{(3)}] a relation set $\mathcal{R}\subset \mathcal{C}\times \mathcal{C}$ having the properties that 
$\mathrm{(a)}$ $v_{\cdot i}\to v_{\cdot i}\notin\mathcal{R}~\forall~ v_{\cdot i}\in\mathcal{C}$; $\mathrm{(b)}~\forall~v_{\cdot i}\in\mathcal{C},~\exists~v'_{\cdot i}\in\mathcal{C}$ either $v_{\cdot i}\to v'_{\cdot i}\in\mathcal{R}$ or $v'_{\cdot i}\to v_{\cdot i}\in\mathcal{R}$. 
\end{itemize}
The triple $\mathcal{(S,C,R)}$ is usually used to represent a CRN.
\end{definition}

Each relation $v_{\cdot i}\to v'_{\cdot i}$ in the network is called a reaction with $v_{\cdot i}$ termed reactant complex while $v'_{\cdot i}$ termed resultant complex, which is equivalent to 
\begin{equation*}
   \sum_{j=1}^{n} v_{ji}S_{j}\longrightarrow\sum_{j=1}^{n}
v'_{ji}S_{j}.
\end{equation*}
A CRN can be also viewed as a directed graph with vertices to represent complexes while directed edges to correspond to reactions. A connected component of the graph is called a linkage class.  

\begin{definition}[Reversible and Weakly Reversible CRN]\label{def:2reversible} A CRN $(\mathcal{S,C,R})$ is called
\begin{itemize}
  \item reversible if $v_{\cdot i}\to v'_{\cdot i}\in\mathcal{R}$ whenever $v'_{\cdot i}\to v_{\cdot i}\in\mathcal{R}$. 
\item  weakly reversible if for any reaction $v_{\cdot i}\to v'_{\cdot i}\in\mathcal{R}$ there exists a series of reactions starting from $v'_{\cdot i}$ and ending with $v_{\cdot i}$, i.e., $v'_{\cdot i} \to v_{\cdot i_{1}}\in\mathcal{R}$, $\cdots$, $v_{\cdot i_{m}}\to v_{\cdot i}\in\mathcal{R},~m\textless r$.  
\end{itemize}  
 \end{definition}

A reversible CRN must be weakly reversible CRN, but not vice versa. 

\begin{definition}[Dimension of Network] For a CRN $(\mathcal{S,C,R})$, $v'_{\cdot i}-v_{\cdot i}$ is called the reaction vector of the $i$th $(i=1,\cdots,r)$ reaction $v_{\cdot i}\to v'_{\cdot i}$. The rank of the set of reaction vectors is the dimension of this network, denoted by $s$, i.e., $\mathrm{rank}\{v'_{\cdot 1}-v_{\cdot 1},\cdots,v'_{\cdot r}-v_{\cdot r}\}=s$.
\end{definition}
\begin{definition}[Stoichiometric Compatibility Class] For a CRN $(\mathcal{S,C,R})$, the linear space spanned by all reaction vectors is named as the stoichiometric subspace, denoted by $\mathscr{S}=\mathrm{span} \{v'_{\cdot 1}-v_{\cdot 1},\cdots,v'_{\cdot r}-v_{\cdot r}\}$. Let $x_0\in\mathbb{R}^\mathcal{S}_{\geq 0}$, the sets $\mathscr{S}(x_0)\triangleq \{x_0+\xi|\xi\in\mathscr{S}\}$, $\mathscr{S}(x_0)\cap\mathbb{R}_{\geq 0}^{\mathcal{S}}$ and  $\mathscr{S}(x_0)\cap\mathbb{R}_{>0}^{\mathcal{S}}$ are called the stoichiometric compatibility class, non-negative stoichiometric compatibility class and positive stoichiometric compatibility class of $x_0$, respectively.
\end{definition}

The reaction vectors define a stoichiometric matrix $\varGamma=(v'_{\cdot 1}-v_{\cdot 1},\cdots,v'_{\cdot r}-v_{\cdot r})\in\mathbb{Z}^{\mathcal{S}\times\mathcal{R}}$. Clearly, $\mathrm{dim}~\mathscr{S}=\mathrm{dim}(\mathrm{Im}~\varGamma)=s$. Geometrically, a CRN may be represented through the Cartesian coordinate system, where a complex is projected to a vector.

\newtheorem{example}{Example}
\begin{example}\label{ex:1}
Consider the following network with two species and four reactions
\begin{align*}
    S_{1}& \leftrightharpoons S_{2},\\
    S_{1}+S_{2}&\to 2S_{2}\to  2S_{1}.
\end{align*}
It is easy to write
$$\mathcal{S}=\{S_{1},S_{2}\},~\mathcal{C}=\left\{\begin{bmatrix}
1 \\0
\end{bmatrix},\begin{bmatrix}
0\\1
\end{bmatrix},\begin{bmatrix}
1\\1
\end{bmatrix},\begin{bmatrix}
0\\2
\end{bmatrix},
\begin{bmatrix}
2\\0
\end{bmatrix}\right\},$$
$$\mathcal{R}=\left\{\begin{bmatrix}
1 \\0
\end{bmatrix}\to\begin{bmatrix}
0\\1
\end{bmatrix},\begin{bmatrix}
0 \\1
\end{bmatrix}\to\begin{bmatrix}
1\\0
\end{bmatrix},\begin{bmatrix}
1\\1
\end{bmatrix}\to\begin{bmatrix}
0\\2
\end{bmatrix},\begin{bmatrix}
0 \\2
\end{bmatrix}\to\begin{bmatrix}
2\\0
\end{bmatrix}\right\},$$
$\mathscr{S}=\mathrm{span}\{(-1,1)^\top\}$ and  $\mathrm{dim}~\mathscr{S}=1$. Therefore, the network is a $1d$ (1-dimensional) CRN and contains two linkage classes. The geometric representation of this network is shown in \textbf{Figure \ref{fig:2.1}}. Clearly, the stoichiometric subspace is a straight line and all reaction vectors are parallel to this line.
\end{example}

Utilizing the geometric property of CRNs, Craciun et al. \cite{Craciun2013} defined the concept of endotactic CRN, which serves closely for the purpose of addressing the persistence issue. They also stated that the class of endotactic networks is larger than the well-known class of weakly reversible networks. Instead, Gopalkrishnan et al. \cite{G2014} redefined this concept using algebraic language and further proposed the notion of strongly endotactic CRN. Their redefinition begins with defining a partial order relation in $\mathbb{R}^{\mathcal{S}}$.    

 \begin{definition}[$w$-Preorder \cite{G2014}] 
Let $w$ be a vector in $\mathbb{R}^{\mathcal{S}}$. The $w$-Preorder on $\mathbb{R}^{n}$, denoted by $\leq_{w}$, is defined by 
 \begin{equation*}
        y_{1} \leq_{w} y_{2}~~\Longleftrightarrow~~\langle w,~y_{1}\rangle \leq \langle w,~y_{2}\rangle,
\end{equation*} 
where $y_1,y_2\in\mathbb{R}^{\mathcal{S}}$, and $\langle \cdot,\cdot\rangle$ represents standard inner product. We write $y_{1}<_{w}y_{2}$ if $\langle w,y_{1}\rangle <\langle w,y_{2}\rangle.$
\end{definition}

\begin{definition}[$\leq_{w}$-$maximal$ and $\leq_{w}$-$minimal$ Elements \cite{G2014}]
Let $w\in\mathbb{R}^{\mathcal{S}}$ and $Y \subset \mathbb{R}^{\mathcal{S}}$. The element $y_r\in Y$ is said to be a $\leq_{w}$-$maximal$ element in $Y$ if 
\begin{equation*}
      \langle w, y_r \rangle \geq \langle w,y\rangle , ~\forall~ y \in Y.
   \end{equation*}
The element $y_l\in Y$ is called a $\leq_{w}$-$minimal$ element in $Y$ if
\begin{equation*}
      \langle w, y_l \rangle \leq \langle w,y\rangle , ~\forall~ y \in Y.
   \end{equation*}    
\end{definition}

When applying this definition to a CRN $(\mathcal{S,C,R})$, a complex is called leftmost relative to $w$ if it is $\leq_{w}$-$minimal$ in $\mathcal{C}$, and rightmost relative to $w$ if it is $\leq_{w}$-$maximal$ in $\mathcal{C}$. We then give the definition of endotactic CRN in $1D$ case, where $1D$ means that the network is a $1d$ CRN with single linkage class or empty set. Geometrically, all complexes in a $1D$ network belong to the same line \cite{Craciun2013}.  

\begin{definition}[$1D$ Endotactic Network \cite{Craciun2013,G2014}]\label{def:1endotatic}
A $1D$ network $(\mathcal{S,C,R})$ is endotactic if 
\begin{itemize}
  \item[\rm{(1)}] it contains empty complex; 
\end{itemize}  
or
\begin{itemize}
  \item[\rm{(2)}] it contains at least two reactant complexes. Moreover, each reaction with a leftmost reactant pointing to the right while each reaction with a rightmost reactant pointing to the left. Here, the preorder $\leq_{w}$ in $\mathbb{R}^{\mathcal{S}}$ is defined by a basis $w\in\mathbb{R}^{\mathcal{S}}_{\neq 0}$ of the stoichiometric subspace $\mathscr{S}$ of the network.
\end{itemize}  
\end{definition}

\begin{example}\label{ex:2gao}
The sub-network of \ref{ex:1}
\begin{align*}
    S_{1}+S_{2}&\to 2S_{2}\to  2S_{1}
\end{align*}
is a $1D$ endotactic CRN since the two extreme source complexes (leftmost or rightmost reactant), $S_1+S_2$ and $2S_2$, both react towards direction of the other source complex, as can be seen in \ref{fig:2.1}.  
 \begin{figure}
        \centering
        \begin{minipage}[c]{0.45\textwidth}
        \centering
        \includegraphics[width=0.8\textwidth]{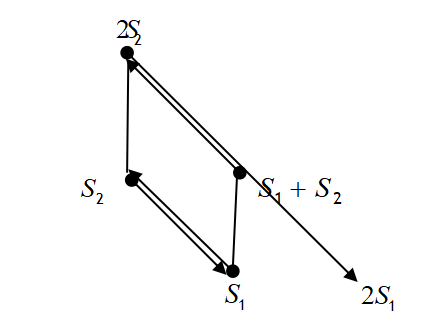}
        \end{minipage}
        \begin{minipage}[c]{0.45\textwidth}
        \centering
          \includegraphics[width=0.8\textwidth]{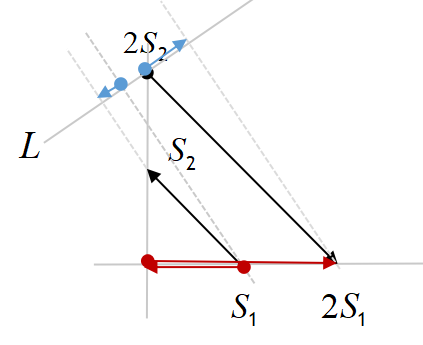}
        \end{minipage}
         \caption{The left graph shows the complexes, reactions and the reactant complexes polytope of \textbf{Example \ref{ex:1}} network in the phase space. The right figure is the first network in \textbf{Example \ref{ex:2}} projecting into the straight line $L$.}
        \label{fig:2.1}
    \end{figure}
\end{example}

\begin{definition}[Endotactic CRN \cite{Craciun2013,G2014}]\label{def:Endotatic} A CRN $\mathcal{(S,C,R)}$ is
 \begin{enumerate}
    \item [\rm{(1)}] $w$-endotactic with respect to a certain $w\in \mathbb{R}^{\mathcal{S}}$ if for any reactant $v_i$ that is $\leq_{w}$-$maximal$ among all reactants such that the reaction vectors are not orthogonal to $w$, we have 
     \begin{equation*}
     \langle w, v'_{\cdot i}-v_{\cdot i}\rangle<0.
 \end{equation*}
 \item [\rm{(2)}] an endotactic network if it is $w$-endotactic with respect to any $w\in \mathbb{R}^{\mathcal{S}}$. 
\item [\rm{(3)}] strongly endotactic if it is endotactic and for each $w$ not orthogonal to its stoichiometric subspace $\mathscr{S}$, there exists a reaction $v_{\cdot i}\to v'_{\cdot i}$ such that:
\begin{itemize}
\item $\langle w, v'_{\cdot i}-v_{\cdot i}\rangle <0$ and
\item $v_{\cdot i}$ is $\leq_{w}$-$maximal$ among all reactants. 
\end{itemize}
\end{enumerate}
\end{definition}
It is clear that among the above three kinds of networks, the set of $w$-endotactic networks is largest, and the set of strongly endotactic networks is smallest. Geometrically, a network is $w$-endotactic if the projection of this network to the line which $w$ belongs to is endotactic, and a network is endotactic if its projection on any line is endotactic \cite{Craciun2013}. The strongly endotactic network means that there exists a nontrival reaction for both leftmost and rightmost reactants in the projection of the original network to any line not orthogonal to its stoichiometric subspace \cite{G2014}.   

\begin{example}\label{ex:2}
The following CRN
\begin{align*}
S_{1}&\to S_{2},\\
2S_{2}&\to 2S_{1}
\end{align*}
is also a $1d$ network with $2$ linkage classes. By projecting the network to the horizontal line containing $w=(1,0)^\top$, as can be seen in \textbf{Fig. \ref{fig:2.1}}, the leftmost reactant complex $2S_{2}$ points to right and the rightmost reactant complex $S_{1}$ points to left. Hence, this network is at least a $w$-endotactic network. We further project it into the line $L$ containing another $w$, and the projection is shown not endotactic.

A look back at \textbf{Example \ref{ex:1}} might reveal that whatever $w$ is, the projection of the network to the line containing $w$ is endotactic, so the network given in \textbf{Example \ref{ex:1}} is endotactic. We also observe \textbf{Example \ref{ex:2gao}} again. It is obviously an endotatic network, and moreover,  for any line $L$ not orthogonal to its stoichiometric subspace, both the leftmost and the rightmost reactants of the projected $1D$ network on $L$ contain nontrival reactions. Therefore, \textbf{Example \ref{ex:2gao}} is a strongly endotactic network.
\end{example}

In fact, it is not easy to check a network endotactic or not due to  arbitrariness of $w$, which also applies to check a strongly endotactic network. Craciun et. al. \cite{Craciun2013} presented the sufficient and necessary conditions for $2$-species endotactic network while Johnson et. al. \cite{M2016} developed a computation algorithm for endotactic and strongly endotactic networks checking.
  
The endotactic structure has a close connection with persistence, a concept from population dynamics. The dynamics of a CRN system capturing the change of concentration of every species $S_j~ (j=1,\cdots,n)$, identified by $x_j$, is given once the reaction rate is specified as a function $V:\mathbb{R}^n_{\geq 0}\mapsto\mathbb{R}^r_{\geq 0}$ of $x=(x_1,\cdots,x_n)^\top$. The most frequently-used model to specify the reaction rate is mass-action kinetics, under which the reaction rate follows power law with respect to the concentration of every species in the reactant complex. For example, for the $i$th reaction $v_{\cdot i}\to v'_{\cdot i}$, the reaction rate $V_i(x)$ is evaluated by
\begin{equation}
 V_i(x)=k_ix^{v_{\cdot i}}\triangleq k_i\prod_{j=1}^n x_j^{v_{ji}}, 
\end{equation}
where $k_i\in\mathbb{R}_{> 0}$ is the reaction rate coefficient. 

\begin{definition}[Mass-action System]
A CRN $\mathcal{(S,C,R)}$ equipped with mass-action kinetics is called a mass-action system (MAS), which is usually represented by the quaternary $(\mathcal{S,C,R},k)$. 
\end{definition}

The dynamical equation of a MAS $(\mathcal{S,C,R},k)$ can be written as
\begin{equation}\label{eq:de1}   
\frac{dx(t)}{dt}\triangleq\dot{x}(t)=\sum_{i=1}^{r}k_{i}x^{v_{\cdot i}}(v'_{\cdot i}-v_{\cdot i}).
\end{equation}

\begin{definition}[Equilibrium]
For a MAS $(\mathcal{S,C,R},k)$, a concentration vector $x^*\in\mathbb{R}^{n}_{>0}$ is called an equilibrium if
    \begin{equation*}
\sum_{i=1}^{r}k_{i}(x^*){^{v_{\cdot i}}(v'_{\cdot i}-v_{\cdot i})}=\mathbbold{0}_n,
\end{equation*}
and a complex balanced equilibrium if
   \begin{equation*}
    \sum_{\{v_{\cdot i}=C\}}k_{i}(x^*)^{v_{\cdot i}}=\sum_{\{v'_{\cdot i}=C\}}k_{i}(x^*)^{v_{\cdot i}},~~~~\forall ~C\in\mathcal{C}.
\end{equation*} 
\end{definition}

A complex balanced equilibrium must be an equilibrium, but not vice verse. A MAS admits an (complex balanced) equilibrium is called (complex) balanced MAS. If there exists a complex balanced equilibrium in a MAS, any other equilibrium (if exists) in this MAS is also a complex balanced equilibrium \cite{Horn1972}.

\section{Persistence}\label{sec:3}
In this section, the concept of persistence and some known results to suggest persistence are introduced \cite{Feinberg1987,Angeli2007,Anderson2008,Anderson2010,Pantea2012,Craciun2013,G2014}. We also introduce the Lyapunov Function PDEs \cite{Fang} that will be used to analyze persistence subsequently.  
\subsection{Definitions and Related Results}
In CRNT, persistence means that none of the concentrations of species can tend to zero if they are not zero at the beginning of reactions. Mathematically, it is defined as follows. 

\begin{definition}[Persistence]
Consider a MAS $(\mathcal{S,C,R},k)$ described by {(\ref{eq:de1})}. This MAS is called persistent if any forward trajectory $x(t)\in \mathbb{R}^{\mathcal{S}}_{\geq 0}$ with positive initial condition $x(0)\in \mathbb{R}^{\mathcal{S}}_{> 0}$ satisfies
\begin{equation*}
    \liminf_{t\rightarrow \infty}x_{j}(t)>0~~~for~all~j\in\{1,\cdots,n\}.
\end{equation*}
\end{definition}
The above definition works for all cases, including the cases of bounded trajectory and of unbounded trajectory. In the case of bounded trajectory, the definition may be reduced to describe its $\omega$-$limit$ point. 

\begin{definition}[$\omega$-$limit$ Point]
The set of $\omega$-$limit$ points for the trajectory $x(t)$ with positive initial condition $x(0)\in \mathbb{R}^{\mathcal{S}}_{>0}$ is
\begin{equation}
    \omega(x(0)):=\{x\in \mathbb{R}^{\mathcal{S}}_{\geq 0}~|~x(t_{N})\rightarrow x,~\text{for some sequence} ~t_{N}\rightarrow \infty\}.
\end{equation}
\end{definition}
\begin{definition}[Persistence for Bounded Trajectory]
For a MAS $(\mathcal{S,C,R},k)$ with bounded trajectory, it is persistent if
\begin{equation}
 \omega(x(0))\cap \partial \mathbb{R}^{\mathcal{S}}_{\geq 0}=\varnothing, ~~~~\forall~x(0)\in \mathbb{R}^{\mathcal{S}}_{>0}. 
\end{equation} 
\end{definition}

The following concepts play an important role on characterizing persistence of a MAS.

\begin{definition}[Semilocking set and Locking Set \cite{Anderson2010,Angeli2007}]\label{def:semilocking}
Consider a CRN $(\mathcal{S,C,R})$. A nonempty set $W\subseteq\mathcal{S}$ is said to be a semilocking set or siphon if one species of $W$ is in a resultant complex $v'_{\cdot i}$, there must exist a species of $W$ in the reactant complex $v_{\cdot i}$. If for any reaction $v_{\cdot i}\to v'_{\cdot i}$ there exists a $S_{j}\in W$ in reactant complex $v_{\cdot i}$, then $W\subseteq \mathcal{S}$ is called a locking set.  
\end{definition}

\begin{definition}[Face \cite{Anderson2010}]\label{def:face}
Consider a CRN $(\mathcal{S,C,R})$ and a nonempty index set $W\subseteq  \mathcal{S}$. The set $Z_{W}$, defined by
\begin{equation*}
    Z_{W}=\{x\in \mathbb{R}^{\mathcal{S}}_{\geq 0}|x_{i}=0~for~S_i\in W \},
\end{equation*}
is named the face of $\mathbb{R}^{\mathcal{S}}_{\geq 0}$. The set $L_W$ given by
\begin{equation*}
    L_W=\{x\in \mathbb{R}^{\mathcal{S}}_{\geq 0}|x_{i}=0~for~S_i\in W~and~x_i>0~for~S_i\notin W\}
\end{equation*}
is called the relative interior of $Z_W$.
\end{definition}

As far as the stoichiometric compatibility class in $(\mathcal{S,C,R})$ is concerned, its face under a given $W$ is
$F_{W}=\mathscr{S}(x)\cap Z_{W},$
and the relative interior of $F_W$ is $\mathscr{S}(x)\cap L_{W}$. Angeli et al. \cite{Angeli2007} also reported a necessary and sufficient condition to suggest a semilocking set. That is:

\begin{proposition}
A nonempty set $W \subseteq \mathcal{S}$ in a CRN $(\mathcal{S,C,R})$ is a semilocking set if and only if the face $F_{W}$ is forward invariable for the dynamics {(\ref{eq:de1})}. 
\end{proposition}
 
The physical explanation of this proposition is easily caught as the concentrations of species in $W$ are always zero and every reaction induced by the species in $W$ will stop.

Based on the concepts of semilocking set, face and endotactic network, many checkable criteria have been reported for persistence, as stated in the Introduction. In this paper, we continue to follow this project, but use another strategy of Lyapunov Function PDEs \cite{Fang}.

\subsection{Lyapunov Function PDEs}
The Lyapunov Function PDEs are initially proposed aiming at producing a Lyapuonov function for capturing the asymptotic stability of equilibria in MASs. This technique is activated by the work of Anderson et. al. \cite{Anderson2015} in which they obtained the  pseudo-Helmholtz free energy function by evaluating an appropriate limit of the stationary distributions for the stochastic MASs. Zhou and Gao \cite{Fang} took full advantage of this evaluation, and approximated the scaling non-equilibrium potential by a unknown function, which actually acts as a potential Lyapunov function. By bridging between the microscopic and the macroscopic level, thermodynamics and potential theory, the Lyapunov function PDEs are yielded from Chemical Master Equation \cite{Fang}. For a MAS $(\mathcal{S,C,R},k)$ governed by {(\ref{eq:de1})}, they include a first-order PDE  
\begin{equation}\label{eq:pde}
    \sum_{i=1}^{r}k_{i}x^{v_{\cdot i}}-\sum_{i=1}^{r}k_{i}x^{v_{\cdot i}}\exp\{(v'_{\cdot i}-v_{.i})^{\top}\triangledown f(x)\}=0,~ x\in \mathbb{R}^{\mathcal{S}}_{>0}
\end{equation}
and a corresponding boundary condition 
\begin{equation}\label{eq:pdea}
   \lim\limits_{x\rightarrow \bar{x}} \sum_{\mathrm{supp}~v_{\cdot i}\subseteq \mathrm{supp}~\bar{x}}k_{i}x^{v_{\cdot i}}-\lim\limits_{x\rightarrow\bar{x}}\sum_{\mathrm{supp}~v'_{\cdot i}\subseteq \mathrm{supp}~\bar{x}}k_{i}x^{v_{\cdot i}}\exp\{(v'_{\cdot i}-v_{\cdot i})^{\top}\triangledown f(x)\}=0,~x\in \mathbb{R}^{\mathcal{S}}_{>0},
\end{equation}
where $\bar{x}\in\mathbb{R}^\mathcal{S}_{\geq 0}$ represents any boundary point and $x-\bar{x}\in\mathscr{S}$.

\begin{remark}
The following boundary condition is straightforward if {(\ref{eq:pdea})} is subtracted from {(\ref{eq:pde})}, i.e.,
\begin{equation}\label{eq:pdeb}
   \lim\limits_{x \rightarrow \bar{x}} \sum_{\mathrm{supp}~v_{\cdot i}\nsubseteq \mathrm{supp}~\bar{x}}k_{i}x^{v_{\cdot i}}-\lim\limits_{x\rightarrow\bar{x}}\sum_{\mathrm{supp}~v'_{\cdot i}\nsubseteq \mathrm{supp}~\bar{x}}k_{i}x^{v_{\cdot i}}\exp\{(v'_{\cdot i}-v_{\cdot i})^{\top}\triangledown f(x)\}=0,~x\in \mathbb{R}^{\mathcal{S}}_{>0}.
\end{equation}
For distinction, {(\ref{eq:pdea})} and {(\ref{eq:pdeb})} are referred to as the boundary condition (a) and (b) in the context, respectively.
\end{remark}

There is a good property about the solutions (if exist) of the PDEs {(\ref{eq:pde}) and (\ref{eq:pdea})}. We exhibit it through the following proposition.

\begin{proposition}[Dissipation of Solutions of the PDEs \cite{Fang}]\label{pro:dissipation}
For a MAS $(\mathcal{S,C,R},k)$ described by {(\ref{eq:de1})}, let $f\in \mathscr{C}^{1}(\mathbb{R}^{\mathcal{S}}_{>0};\mathbb{R})$ be a solution of its Lyapunov Function PDEs {(\ref{eq:pde}) and (\ref{eq:pdea})}, then 
 \begin{equation}
     \frac{df(x)}{dt}=\dot{f}(x)\leq 0,~~x\in \mathbb{R}^{\mathcal{S}}_{>0}
 \end{equation}
with equality to hold if and only if $\triangledown f(x)\perp \mathscr{S}$. Further, assume that $f\in \mathscr{C}^{2}(\mathbb{R}^{\mathcal{S}}_{>0};\mathbb{R})$ is a solution of the PDEs. If $\exists~ D\subset \mathbb{R}^{\mathcal{S}}_{>0}$ such that $\forall x\in D$ and $\forall \mu \in \mathscr{S}$ there is 
 \begin{equation}
\mu^{\top}\triangledown^{2}f(x)\mu\geq 0~~with~equality~to~ hold ~if~ and ~ only~ if ~\mu=\mathbbold{0}_{n},
 \end{equation}
 then for all $x\in D$, $\dot{f}(x)=0$ if and only if $x$ is an equilibrium of the MAS.
\end{proposition}

Based on the dissipation of solutions of the PDEs, Fang and Gao \cite{Fang} succeed in proving that the PDEs {(\ref{eq:pde}) and (\ref{eq:pdea})} work well for capturing the asymptotic stability of equilibria in complex balanced MASs, all $1d$ MASs and some special MASs with dimension beyond $1$. In this paper, we will use the PDEs to analyze the persistence of MASs.

\section{Lyapunov Function PDEs for persistence}\label{sec:4}
In this section, we will present some criteria for judging persistence based on different cases of solutions of the Lyapunov Function PDEs, which are further applied to analyze different boundary points in MASs. 
\subsection{Criteria for judging persistence}
The Lyapunov Function PDEs {(\ref{eq:pde}) and (\ref{eq:pdea})} might serve for analyzing the persistence of bounded MASs by judging that any boundary point is not a $\omega$-$limit$ point according to different values-taking of solutions (if exist). To this task, we firstly reach a basic but very central theorem about the residence time of trajectories of a MAS in the positive quadrant. 
\begin{theorem}\label{thm:inf}
For a MAS $(\mathcal{S,C,R},k)$ given by {(\ref{eq:de1})}, each bounded trajectory starting from a positive initial point $x(0)\in\mathbb{R}^{\mathcal{S}}_{> 0}$ stays in the positive stoichiometric class of $x(0)$ for infinite time, i.e.,
\begin{equation}\label{eq:ep}
    T_+=\lim\limits_{t\rightarrow +\infty}\int_{0}^{t}\mathbbold{1}_{\{x(s)\in \mathscr{S}(x(0))\cap\mathbb{R}^{\mathcal{S}}_{> 0}\}}ds=+\infty,
\end{equation}
where $\mathbbold{1}_{\{\cdot\}}$ is the indicator function defined by $$\mathbbold{1}_{\{A\}}=
\begin{cases}
          1,~A~\mathrm{is~ true};\\
          0,~A~\mathrm{is~ false}.
   \end{cases}
   $$
\end{theorem}
\noindent{\textbf{Proof.}}
For a MAS $(\mathcal{S,C,R},k)$ given by {(\ref{eq:de1})}, its trajectory $x(t)\in\mathbb{R}^{\mathcal{S}}_{\geq 0}$ appears either in the interior or in the boundary of $\mathscr{S}(x(0))\cap\mathbb{R}^{\mathcal{S}}_{\geq 0}$. There are two kinds of possible boundary points, denoted by $\bar{x}\in\partial_{\mathscr{S}(x(0))\cap\mathbb{R}^{\mathcal{S}}_{\geq 0}}$, according to $W\triangleq \mathrm{supp}^{c}\bar{x}$ being a semilocking set or not. We continue the proof in two separate cases. 

    (1) Case I: $W$ is not a semilocking set. When $x(t)$ enters into the boundary $L_{W}$, defined in {(\ref{def:face})}, we have $x_j(t)=0$, $\forall j\in W$. Since $W$ is not a semilocking set, $\exists~ p\in W$ such that $S_p\in \mathrm{supp}~v'_{\cdot i}$ and $\mathrm{supp}~v_{\cdot i} \subseteq W^{c}$, i.e., for the reaction $v_{\cdot i}\to v'_{\cdot i}$ at time $t$, the inflow of species $S_p$ is larger than zero while its outlow is zero, which means that $\dot{x}_p(t)> 0$ and $x_p(t)=0$. This suggests that the trajectory $x(t)$ will leave $L_{W}$ once it enters into it. Therefore, the time set, given by $\{t|x(t)\in L_{W}\}$, is at most discrete. As time $t$ goes to infinity, the trajectory $x(t)$ will stay out of $L_W$ for infinite time.
    
   (2) Case II: $W$ is a semilocking set. In this case, assume that at time $T$ ($T<+\infty$) the trajectory $x(t)$ enters into $L_W$, i.e., $x(T)=\bar{x}\in L_W$. For simplicity of symbols, denote $\dot{x}(t)=g(x)$ for the dynamics {(\ref{eq:de1})}. This ODE together with the initial condition $x(t)=x(0)$ at $t=0$ has the same solution with the following initial value problem  
 \begin{equation}\label{ODE1}
 \left\{
 \begin{array}{ll}
          \frac{dx(T-t)}{dt}=-g(x(T-t)),~~t\in [0,T];\\
           x(T)=\bar{x}, g(x(T))=g(\bar{x}),~~t=0.
          	\end{array}
          		\right.
   \end{equation} 
Since $g(x(T-t))$ is a polynomial function and $x(T-t)$ is bounded, which support the conditions for existence and uniqueness of solution of the above ODE, the solution written by $x(T-t)=(x_{1}(T-t),\cdots,x_{n}(T-t))^\top$ is unique. Further, let $X_{1}\in\mathbb{R}^{W}_{\geq 0}$ and $X_{2}\in\mathbb{R}^{W^c}_{\geq 0}$, then each $g_{j}(x(T-t))=g_{j}(X_{1}(T-t),X_{2}(T-t)),~j=1,\cdots,n,$ is a linear combination of $k_{1}X_{1}^{v_{\cdot 1}|_{W}}X_{2}^{v_{\cdot 1}|_{W^{c}}},\cdots,k_{r}X_{1}^{v_{\cdot r}|_{W}}X_{2}^{v_{\cdot r}|_{W^{c}}}$, where $v_{\cdot 1}|_{W}$ represents a sub-vector of $v_{\cdot 1}$ with every entry orientated by $W$. By denoting $k_{i}X_{2}^{v_{\cdot i}|_{W^{c}}}$ as $k_{i}(t),~i=1,\cdots,r,$ we rewrite $g_{j}(X_{1}(T-t),X_{2}(T-t))=g_{j}(X_{1}(T-t),t)$. The initial value problem {(\ref{ODE1})} with respect to $X_1$ can be rewritten as  
 $$\begin{cases}
        &\frac{dX_{1}(T-t)}{dt}=-g(X_{1}(T-t),t)|_{W},~~t\in [0,T];\\
        & X_{1}(T)=0, g_{j}(0,0)=0,\forall j\in W,~~ t=0.
 \end{cases}$$
Note that this ODE is sole since the uniqueness of $x(T-t)$ implies $X_1(T-t)$ and $k_i(t)$ to be unique, so is $g(X_{1}(T-t),t)|_{W}$. It is not difficult to verify $g_{j}(0,t)=0$ for all $j\in W$ from the fact that every $g_j(X_{1}(T-t),t)$ is a linear combination of $k_{1}(t)X_{1}^{v_{\cdot 1}|_{W}},\cdots,k_{r}(t)X_{1}^{v_{\cdot r}|_{W}}$. Whatever the species in a semilocking set $W$ acts as a reactant or a resultant, the term $X_1^{v_{\cdot i}|_{W}}=0 ~(i=1,\cdots,r)$. We thus have that $X_{1}(T-t)\equiv 0$ is the unique solution of the above equation. Naturally, at $t=T$ one should get $X_1(0)=0$. This is contradicted against $x(0)\in\mathbb{R}^\mathcal{S}_{> 0}$. Therefore, the assumption that the trajectory $x(t)$ enters into $L_W$ in a finite amount of time is not true. 

By combining (1) and (2), we conclude that each bounded trajectory in a MAS with a positive initial point will stay in the positive stoichiometric compatibility class for infinity time when time trends to infinity.$\Box$ 

The proof process also indicates the time that a bounded trajectory of a MAS stays in the boundary.  
\begin{corollary}\label{cor:inf}
Given a MAS $(\mathcal{S,C,R},k)$ with the dynamics of {(\ref{eq:de1})}, the time set for any bounded trajectory starting from a positive initial point to stay in the boundary $L_W$, given in \textbf{Definition \ref{def:face}}, is at most discrete if $W$ is not a semilocking set. Also, it is impossible for the mentioned trajectory to enter into a semilocking boundary in a finite amount of time.   
\end{corollary} 

\textbf{Theorem \ref{thm:inf}} and \textbf{Corollary \ref{cor:inf}} report definitely the time allocation for a bounded trajectory of a MAS to stay in the positive quadrant and the non-semilocking boundary, respectively. The infinite time for the former implies great possibility for a MAS to be persistent while the discrete time set for the latter might suggest a non-semilocking boundary to be not a convergence region. The difficulty in proving persistence may lie in analyzing semilocking boundary, which has been confirmed in Anderson's work \cite{Anderson2008}. The results of \textbf{Theorem \ref{thm:inf}} and \textbf{Corollary \ref{cor:inf}} keep it unknown whether a bounded trajectory can enter into a semilocking boundary in infinite time. Despite these facts, it is still possible to infer persistence of a MAS through different values-taking of solution of its Lyapunov Function PDEs. 

\begin{lemma}\label{lm:4.1}
Consider a MAS $(\mathcal{S,C,R},k)$ described by {(\ref{eq:de1})}. Assume that $f\in \mathscr{C}^{1}(\mathbb{R}^{\mathcal{S}}_{>0};\mathbb{R})$ is a nonnegative solution of its Lyapunov Function PDEs (\ref{eq:pde}) and (\ref{eq:pdea}), and $\bar{x}\in\mathbb{R}^\mathcal{S}_{\geq 0}$ represents any boundary point. $\bar{x}$ is not a $\omega$-$limit$ point if the solution $f(x)$ satisfies
$$\varlimsup\limits_{x\rightarrow \bar{x}}\dot{f}(x)=-\infty~or~\varlimsup\limits_{x\rightarrow\bar{x}}\dot{f}(x)\leq-M,$$ where $x\in \mathbb{R}^\mathcal{S}_{>0}\cap \mathscr{S}(\bar{x})$ and $M$ is a positive constant.
\end{lemma} 
\noindent{\textbf{Proof.}}
For any initial point $x(0)\in \mathbb{R}^\mathcal{S}_{>0}\cap \mathscr{S}(\bar{x})$, we have $f(x(0))\geq 0$ and  
	\begin{equation}\label{eq:3.6}
		f(x(t))=f(x(0))+\int_{0}^{t} \dot{f}(x(s))ds.
		\end{equation}
In the following, we use rebuttals of evidence to continue the proof.

Assume that $\bar{x}$ is a $\omega$-$limit$ point under the given initial conditions. Under this assumption and the conditions of $f\in \mathscr{C}^{1}(\mathbb{R}^{\mathcal{S}}_{>0};\mathbb{R})$ and $\varlimsup\limits_{x\rightarrow \bar{x}}\dot{f}(x)=-\infty$ (or $-M$), we can find a positive constant $M_0$ and a $\epsilon$-neighbourhood of $\bar{x}$ ($\epsilon>0$), denoted by $N_{\epsilon}(\bar{x})$, for any trajectory $x(t)$ starting from $x(0)$ such that $\dot{f}(x)<-M_0,\forall x \in N_{\epsilon}(\bar{x})$, and its residence time in $N_{\epsilon}(\bar{x})$ is infinite, i.e., 
\begin{equation}\label{eq:rt}
    T_{\epsilon}=\lim\limits_{t\rightarrow +\infty}\int_{0}^{t}\mathbbold{1}_{\{x(s)\in N_{\epsilon}(\bar{x})\}}ds=+\infty.
\end{equation}
If we cannot find such a $N_{\epsilon}(\bar{x})$ supporting {(\ref{eq:rt})}, consider another $\epsilon'$-neighbourhood of $\bar{x}$ ($\epsilon'>\epsilon$), denoted by $N_{\epsilon'}(\bar{x})$, in which $\varlimsup\limits_{x\rightarrow \bar{x}}\dot{f}(x)<-M'_0 ~(M'_0>0)$, and denote the residence time for the trajectory $x(t)$ to stay in $N_{\epsilon'}(\bar{x})$ by
\begin{equation*}
    T_{\epsilon'}=\lim\limits_{t\rightarrow +\infty}\int_{0}^{t}\mathbbold{1}_{\{x(s)\in N_{\epsilon'}(\bar{x})\}}ds.
\end{equation*}
If $T_{\epsilon'}=+\infty$, then replace $\epsilon$ in {(\ref{eq:rt})} by $\epsilon'$. If $T_{\epsilon'}$ is finite, then from \textbf{Theorem \ref{thm:inf}} saying that the trajectory $x(t)$ stays in the positive stoichiometric compatibility class of $x(0)$ for infinite time, we get that the time for $x(t)$ to stay out of $N_{\epsilon'}(\bar{x})$ is infinite. Define an area $\Omega$ by   
$$\Omega\triangleq\{x\in \mathbb{R}^{n}_{>0}|\epsilon\leq\|x-\bar{x}\|\leq \epsilon'\},$$ 
which is obviously a bounded region. Since $\bar{x}$ is assumed to be a $\omega$-limit point of the trajectory $x(t)$, for any $t_1>0$ we can find $t_2>t_1$ such that $x(t_{2})\in N_{\epsilon}(\bar{x})$. This together with the trajectory staying out of $N_{\epsilon'}(\bar{x})$ for infinite time suggests that the trajectory will pass through the region $\Omega$ infinite times. Note that $g(x)$ is a continuous function about $x$ and the trajectory $x(t)$ is bounded, so $g(x)$ is bounded and $\exists~M_1>0$ such that $\Vert g(x)\Vert_{\infty}\leq M_{1},\forall x\in \Omega$. Denote the time for the trajectory to enter into $\Omega$ at a certain time by $t_a$ and to leave $\Omega$ by $t_b$, then we have    
	$$\epsilon'-\epsilon\leq \|x(t_b)-x(t_a)\|=\|\int_{t_a}^{t_b}g(x(s))ds\|\leq(t_b-t_a)M_{1}\Rightarrow t_b-t_a\geq \frac{\epsilon'-\epsilon}{M_{1}},$$
where $g(x(t))=\dot{x}(t)$. This means that it at least needs the time of $(\epsilon'-\epsilon_{1})/M_{1}$ for the trajectory to stay in $\Omega$. Since the trajectory passes through $\Omega$ infinite times, so its residence time in $\Omega$ will be infinite if the time trends to infinity. This is inconsistent with $T_{\epsilon'}$ being finite. 

Therefore, if we assume $\bar{x}\in\omega(x(0))$, there must exist a $N_\epsilon(\bar{x})$ such that $\varlimsup\limits_{x\rightarrow \bar{x}}\dot{f}(x)<-M_0$ and $T_\epsilon=+\infty$. By taking limit on both sides of {(\ref{eq:3.6})}, we get 
\begin{align*}
	    \varliminf\limits_{t\rightarrow +\infty}f(x(t))&=f(x(0))+\lim\limits_{t\rightarrow+\infty}\int_{0}^{t}\dot{f}(x(s))ds\\
	    &\leq f(x(0))-M_{0}\lim\limits_{t\rightarrow +\infty}\int_{0}^{t}\mathbbold{1}_{\{x(s)\in N_{\epsilon}(\bar{x})\}}ds\\&=-\infty,
	\end{align*}
where the second inequality holds since $\dot{f}(x)\leq 0$, as said in \textbf{Proposition \ref{pro:dissipation}}. This result is contrary to the fact that $f(x)\geq 0$. Hence, the assumption $\bar{x}\in\omega(x(0))$ is not true, which completes the proof. $\Box$

\begin{lemma}\label{lm:4.3}
For a MAS $(\mathcal{S,C,R},k)$ given by (\ref{eq:de1}), let $f\in \mathscr{C}^{1}(\mathbb{R}^n_{>0};\mathbb{R})$ be a nonnegative solution of its Lyapunov Function PDEs {(\ref{eq:pde})} and {(\ref{eq:pdea})}, and $\bar{x}\in\mathbb{R}^\mathcal{S}_{\geq 0}$ be any boundary point. If $\bar{x}$ is the local maximal point for the solution $f(x)$, namely, $\lim\limits_{\Gamma\rightarrow \bar{x}}f(x)$ is the maximum value for each path $\Gamma$. Then $\bar{x}$ is not a $\omega$-$limit$ point. 
\end{lemma}
\noindent{\textbf{Proof.}}
We also use rebuttals of evidence to conduct proof. Assume that $\bar{x}$ is a $\omega$-$limit$ point. Namely, for any initial point $x(0)\in \mathbb{R}^\mathcal{S}_{>0}\cap \mathscr{S}(\bar{x})$, there exists a time series $\{t_N\}$ such that $\bar{x}\in\omega(x(0))$. Also, $\forall~\epsilon>0$, $\exists$ a moment $t_\epsilon<t_N\to + \infty$ such that $x(t_\epsilon)\in N_{\epsilon}(\bar{x})$, where $N_{\epsilon}(\bar{x})$ is a $\epsilon$-neighborhood of $\bar{x}$. Note that $\dot{f}(x)\leq 0$, so it is impossible that $f(\bar{x})$ ($f(x)$ can be defined at $x=\bar{x}$) or $\lim\limits_{\substack{x\rightarrow \bar{x}\\x\in x(t)}}f(x)$ ($f(x)$ can not be defined at $x=\bar{x}$) is local maximal. Therefore, $\bar{x}$ is not a $\omega$-$limit$ point. $\Box$

The above results suggest that we can judge a boundary point in a MAS not a $\omega$-$limit$ point utilizing the solutions of its Lyapunov function PDEs. We thus can further deduce persistence of MASs from these results for the bounded trajectories.  

\begin{theorem}\label{Th:4.5}
For a MAS $(\mathcal{S,C,R,},k)$ governed by {(\ref{eq:de1})}, let its Lyapunov function PDEs {(\ref{eq:pde}) and (\ref{eq:pdea})} contain a nonnegative solution $f\in\mathscr{C}^1(\mathbb{R}^n_{>0};\mathbb{R})$, and its trajectories be bounded. Then the MAS is persistent if one of the following conditions holds for any boundary point $\bar{x}\in\mathbb{R}^{\mathcal{S}}_{\geq 0}$:
\begin{itemize}
		\item[\rm{(1)}] $\varlimsup\limits_{x\rightarrow \bar{x}}\dot{f}(x)=-\infty$;
			\item[\rm{(2)}] $\varlimsup\limits_{x\rightarrow \bar{x}}\dot{f}(x)=-M$ with $M$ a positive constant;
			\item[\rm{(3)}] $\bar{x}$ is the local maximal point for $f(x)$.
		\end{itemize}
	\end{theorem}
\noindent{\textbf{Proof.}}
Combining \textbf{Lemma \ref{lm:4.1}} and \textbf{\ref{lm:4.3}}, we have that if one the above three conditions holds, then $\bar{x}$ is not a $\omega$-$limit$ point. Since $\bar{x}$ is any boundary point, the result means that
$\omega(x(0))\cap \partial \mathbb{R}^{\mathcal{S}}_{\geq 0}=\varnothing$ for any initial point $x(0)\in \mathscr{S}(\bar{x})\cap \mathbb{R}^{\mathcal{S}}_{>0}$. This together with the known condition of bounded trajectories supports the MAS to be persistent. $\Box$

When addressing the issue of persistence, a key point is to identify whether the set related to boundary points $\bar{x}$, defined through $W\triangleq \mathrm{supp}^{c}\bar{x}$, is a semilocking set or not. The result is quite clear in the case of $W$ being not a semilocking set. Angeli et al. \cite{Angeli2007} and Anderson \cite{Anderson2008} has proved independently that a MAS is persistent if $W$ is not a semilocking set for any boundary point. In the following, we will revisit this result through the Lyapunov function PDEs method, and discuss some possibilities to reach persistence in the case of $W$ being a semilocking set.

\subsection{For non-semilocking boundary points}
For non-semilocking boundary points $\bar{x}\in\mathbb{R}^{\mathcal{S}}_{\geq 0}$, i.e., $W\triangleq \mathrm{supp}^{c}\bar{x}$ is not a semilocking set, the following lemma characterizes the time derivative of solution of the Lyapunov function PDEs as the trajectory trends to the non-semilocking boundary point.
	\begin{lemma}\label{lm:3.2}
	Consider a MAS $(\mathcal{S,C,R},k)$ with $f\in \mathscr{C}^{1}(\mathbb{R}^n_{>0};\mathbb{R})$ as a solution of its Lyapunov function PDEs (\ref{eq:pde}) and (\ref{eq:pdea}). Let $\bar{x}\in\mathbb{R}^{\mathcal{S}}_{\geq 0}$ represent any boundary point of this MAS. If the set $W\triangleq \mathrm{supp}^{c}\bar{x}$ is not a semilocking set, then $\varlimsup\limits_{x\rightarrow \bar{x}}\dot{f}(x)=-\infty$ where $x\in \mathbb{R}^{n}_{>0}$ is bounded and satisfies $x-\bar{x}\in\mathscr{S}$. 
		\end{lemma}
\noindent{\textbf{Proof.}}
The time derivative of $f(x)$ follows 
\begin{align}\label{fdot}
\dot{f}(x(t)) & =\sum_{i=1}^{r}k_{i}x^{v_{\cdot i}}(v'_{\cdot i}-v_{\cdot i})^{\top}\cdot\triangledown f(x)\\
& =\sum_{\mathrm{supp}~v'_{\cdot i}\subseteq \mathrm{supp}~ \bar{x}}k_{i}x^{v_{\cdot i}}(v'_{\cdot i}-v_{\cdot i})^{\top}\triangledown f(x)+\sum_{\mathrm{supp}~v'_{\cdot i}\nsubseteq \mathrm{supp}~ \bar{x}}k_{i}x^{v_{\cdot i}}(v'_{\cdot
i}-v_{\cdot
i})^{\top}\triangledown f(x).\notag
\end{align}  
Since $f\in \mathscr{C}^{1}(\mathbb{R}^n_{>0},\mathbb{R})$ is a solution of the Lyapunov function PDEs, it should satisfy the boundary conditions {(\ref{eq:pdea}) and (\ref{eq:pdeb})}.

For the boundary condition (a) {(\ref{eq:pdea})}, we have   
\begin{equation*}
{\lim_{x\rightarrow \bar{x}}\sum_{\mathrm{supp}~v_{\cdot i}\subseteq \mathrm{ supp}~\bar{x}}k_{i}x^{v_{\cdot i}}}={\lim_{x\rightarrow \bar{x}}\sum_{\mathrm{supp} ~v'_{\cdot i}\subseteq \mathrm{supp}~ \bar{x}} k_{i}x^{v_{\cdot i}}\exp \{(v'_{\cdot i}-v_{\cdot i})^{\top}\triangledown f(x)\}}.
\end{equation*}
By letting $a=(v'_{\cdot i}-v_{\cdot i})^\top\triangledown f(x)$ and utilizing the fact of $e^{a}\geq a+1$, we get
\begin{equation*}\label{eq:a}
{\varlimsup\limits_{x\rightarrow \bar{x}}(\sum_{\substack{\mathrm{supp}~ v_{\cdot i}\subseteq \mathrm{supp}~\bar{x}\\\mathrm{supp}~ v'_{\cdot i}\nsubseteq \mathrm{supp}~ \bar{x}}}k_{i}x^{v_{\cdot i}}-\sum_{\substack{\mathrm{supp}~v_{\cdot i}\nsubseteq \mathrm{supp}~\bar{x}\\\mathrm{supp}~v'_{\cdot i}\subseteq \mathrm{supp} ~\bar{x}}}k_{i}x^{v_{\cdot i}})}\geq \varlimsup_{x\rightarrow \bar{x}}\sum_{\mathrm{supp}~v'_{\cdot i}\subseteq \mathrm{supp}~\bar{x}}k_{i}x^{v_{\cdot i}}(v'_{\cdot i}-v_{\cdot i})^{\top}\triangledown f(x).
\end{equation*}
Since $W\triangleq \mathrm{supp}^{c}\bar{x}$ is not a semilocking set, there at least exists a reaction $v_{\cdot i}\rightarrow v'_{\cdot i}$ which satisfies $\mathrm{supp}~ v_{\cdot i}\subseteq \mathrm{supp}~\bar{x}$ and $\mathrm{supp}~ v'_{\cdot i}\nsubseteq \mathrm{supp}~\bar{x}$. This together with the boundedness of $x(t)$ suggests that there exists a constant $M_1\geq 0$ supporting
$$\lim \limits_{x\rightarrow \bar{x}}\sum_{\substack{\mathrm{supp} ~v_{\cdot i}\subseteq \mathrm{supp}~\bar{x}\\\mathrm{supp} ~v'_{\cdot i}\subseteq \mathrm{supp}~\bar{x}}}k_{i}x^{v_{\cdot i}}= M_1~\mathrm{while}~\lim_{x\rightarrow \bar{x}}\sum_{\substack{\mathrm{supp}~v_{\cdot i}\nsubseteq \mathrm{supp}~\bar{x}\\\mathrm{supp}~v'_{\cdot i}\subseteq \mathrm{supp} ~\bar{x}}}k_{i}x^{v_{\cdot i}}=0.$$
We thus have 
\begin{equation}\label{eq:3.2}
{\varlimsup\limits _{x\rightarrow \bar{x}}\sum_{\mathrm{supp}~ v'_{\cdot i}\subseteq \mathrm{supp}~\bar{x}}k_{i}x^{v_{\cdot i}}(v'_{\cdot i}-v_{\cdot i})^{\top}\triangledown f(x)\leq M_1.}
\end{equation}

For the boundary condition (b) {(\ref{eq:pdeb})}, we have 
\begin{equation}\label{gaobd}
\lim\limits_{x\rightarrow \bar{x}}\sum_{\mathrm {supp}~v_{\cdot i}\nsubseteq \mathrm{ supp}~\bar{x}}k_{i}x^{v_{\cdot i}}=\lim\limits_{x\rightarrow \bar{x}}\sum_{\mathrm{supp}~v'_{\cdot i}\nsubseteq \mathrm{supp}~\bar{x}}k_{i}x^{v_{\cdot i}}\exp\{(v'_{\cdot i}-v_{\cdot i})^{\top}\triangledown f(x)\}. 
\end{equation}
Note that the left terms are zero while the right terms can be decomposed into \begin{equation*}
\lim\limits_{x\rightarrow \bar{x}}\sum_{\substack{\mathrm{supp}~v_{\cdot i}\nsubseteq \mathrm{supp}~\bar{x}\\\mathrm{supp}~v'_{\cdot i}\nsubseteq \mathrm{supp} ~\bar{x}}}k_{i}x^{v_{\cdot i}}\exp\{(v'_{\cdot i}-v_{\cdot i})^{\top}\triangledown f(x)\}+\lim\limits_{x\rightarrow \bar{x}}\sum_{\substack{\mathrm{supp}~v_{\cdot i}\subseteq \mathrm{supp}~\bar{x}\\\mathrm{supp}~v'_{\cdot i}\nsubseteq \mathrm{supp} ~\bar{x}}}k_{i}x^{v_{\cdot i}}\exp\{(v'_{\cdot i}-v_{\cdot i})^{\top}\triangledown f(x)\},
\end{equation*}  
where both parts are nonnegative, so we get
\begin{equation*}
\lim\limits_{x\rightarrow \bar{x}}\sum_{\substack{\mathrm{supp}~v_{\cdot i}\nsubseteq \mathrm{supp}~\bar{x}\\\mathrm{supp}~v'_{\cdot i}\nsubseteq \mathrm{supp} ~\bar{x}}}k_{i}x^{v_{\cdot i}}\exp\{(v'_{\cdot i}-v_{\cdot i})^{\top}\triangledown f(x)\}=0.
\end{equation*}   
Also, by applying $e^{a}\geq a+1$ with $a=(v'_{\cdot i}-v_{\cdot i})^{\top}\triangledown f(x)$ to the above equation, we obtain 
$$\varlimsup\limits_{x\rightarrow \bar{x}}\sum_{\substack{\mathrm{supp}~v_{\cdot i}\nsubseteq \mathrm{supp}~\bar{x}\\\mathrm{supp}~v'_{\cdot i}\nsubseteq \mathrm{supp} ~\bar{x}}}k_{i}x^{v_{\cdot i}}(v'_{\cdot i}-v_{\cdot i})^{\top}\triangledown f(x)+\lim\limits_{x\rightarrow\bar{x}}\sum\limits_{\substack{\mathrm{supp}~ v_{\cdot i}\nsubseteq \mathrm{supp}~ \bar{x}\\\mathrm{supp} ~v'_{\cdot i}\nsubseteq \mathrm{supp} ~\bar{x}}}k_{i}x^{v_{\cdot i}}\leq 0,$$
which means 
\begin{equation}\label{eq:3.3}
\varlimsup\limits_{x\rightarrow\bar{x}}\sum_{\substack{\mathrm{supp} ~v_{\cdot i}\nsubseteq \mathrm{supp} ~\bar{x}\\\mathrm{supp}~ v'_{\cdot i}\nsubseteq \mathrm{supp} ~\bar{x}}}k_{i}x^{v_{\cdot i}}(v'_{\cdot i}-v_{\cdot i})^{\top}\triangledown f(x)\leq 0.
\end{equation}
Further, from {(\ref{gaobd})} we get that for all reactions:~$v_{\cdot i}\rightarrow v'_{\cdot i}$ satisfying $\mathrm{supp}~v'_{\cdot i}\nsubseteq \mathrm{supp}~\bar{x}$ and $\mathrm{supp}~v_{\cdot i}\subseteq \mathrm{supp}~\bar{x}$ there are
$$k_{i}x^{v_{\cdot i}}\exp\{(v'_{\cdot i}-v_{\cdot i})^{\top}\triangledown f(x)\}\leq \sum_{\mathrm{supp}~v_{\cdot i}\nsubseteq \mathrm{supp}~\bar{x}}k_{i}x^{v_{\cdot i}}\triangleq F(x) $$
and ${\lim\limits_{x\rightarrow \bar{x}}F(x)=0}$. Hence, we have
$${\varlimsup\limits_{x\rightarrow \bar{x}}~k_{i}x^{v_{\cdot i}}(v'_{\cdot i}-v_{\cdot i})^{\top}\triangledown f(x)}\leq {\lim\limits_{x\rightarrow \bar{x}}(\ln F(x)-\ln (k_{i}x^{v_{\cdot i}}))k_{i}x^{v_{\cdot i}}}.$$
Note that $W$ is not a semi-locking set, so $\lim\limits_{x\rightarrow \bar{x}}\ln F(x)=-\infty$ and $\lim\limits_{x\rightarrow \bar{x}} k_{i}x^{v_{\cdot i}}$ trends to a finite number. Therefore, we have
\begin{equation}\label{eq:3.4}
\varlimsup\limits_{x\rightarrow \bar{x}}\sum_{\substack{\mathrm{supp}~ v_{\cdot i}\subseteq \mathrm{supp} ~\bar{x}\\\mathrm{supp} ~v'_{\cdot i}\nsubseteq \mathrm{supp}~ \bar{x}}}k_{i}x^{v_{\cdot i}}(v'_{\cdot i}-v_{\cdot i})^{\top}\cdot\triangledown f(x)=-\infty. 
\end{equation}
Taking limit from both sides of {(\ref{fdot})} as $x\to\bar{x}$, and combining {(\ref{eq:3.2}), (\ref{eq:3.3}) and (\ref{eq:3.4})}, we get $\varlimsup\limits_{x\rightarrow \bar{x}}\dot{f}(x)=-\infty$. $\Box$

Applying \textbf{Theorem \ref{Th:4.5}} to the case of non-semilocking boundary points, we might reasonably conclude persistence for some MASs.   

\begin{theorem}\label{Th:3.4}
    For a MAS $(\mathcal{S,C,R},k)$, assume its Lyapunov Function PDEs (\ref{eq:pde}) and (\ref{eq:pdea}) admit a nonnegative solution $f\in \mathscr{C}^{1}(\mathbb{R}^n_{>0};\mathbb{R})$. Let $\bar{x}\in\mathbb{R}^{\mathcal{S}}_{\geq 0}$ represent any boundary point of the MAS. If $\bar{x}$ is a non-semilocking boundary point, then any bounded trajectory $x(t)\in\mathbb{R}^{\mathcal{S}}_{> 0}$ satisfying $x-\bar{x}\in\mathscr{S}$ is persistent.
	\end{theorem}
\noindent{\textbf{Proof.}}
The result is immediately by applying \textbf{Theorem \ref{Th:4.5}} and \textbf{Lemma \ref{lm:3.2}} to the given conditions here. $\Box$ 		

It should be mentioned again that the result in \textbf{Theorem \ref{Th:3.4}} is not new, which has appeared in the literature \cite{Angeli2007,Anderson2008}. However, we provide an alternative method here, the Lyapunov Function PDEs, to reach it. We exhibit an application of the current method through the following example.

\begin{example}
A MAS takes the reaction route as
	\begin{equation*}
		\begin{array}{c}
	2S_1\xrightarrow{1}S_1+S_2,\\
	2S_2\xrightarrow{1}S_2+S_3,\\
	2S_3\xrightarrow{1}S_3+S_1.
		\end{array}
		\end{equation*}
It is not difficult to check that for this MAS any nonzero boundary point is a non-semilocking boundary point. The origin is excluded since it is not likely to be an $\omega$-$limit$ point due to mass conservation of this network. Zhou and Gao \cite{Fang} has proved that the following function given by
\begin{equation}
	f(x)=2x^{\top}\mathrm{Ln}\left(\frac{x}{x^{*}}\right)-2\mathbbold{1}_{3}^{\top}(x-x^{*}),~~~x\in\mathbb{R}^n_{>0},
	\end{equation}
is a nonnegative solution of the Lyapunov function PDEs of the MAS, where $x^{*}=(1,1,1)^{\top}$ is an interior equilibrium. Attached to the current MAS, we can compute the limit of its time derivative to be
\begin{align*}
    \lim_{x\rightarrow\bar{x}}\dot{f}(x)&=\lim_{x\rightarrow \bar{x}}2\left[(-x_{1}^{2}+x_{3}^{2})\ln{x_{1}}+(x_{1}^{2}-x_{2}^{2})\ln{x_{2}}+(x_{2}^{2}-x_{3}^{2})\ln{x_{3}}\right].
   \end{align*}
Clearly, for any nonzero boundary point $\bar{x}$, if $\bar{x}-x^*\in\mathscr{S}$, we have $\lim_{x\rightarrow\bar{x}}\dot{f}(x)=-\infty$. This means that the MAS has persistence further based on \textbf{Theorem \ref{Th:3.4}}.
\end{example}

\subsection{For semilocking boundary points}
It has been discussed in \textbf{Corollary \ref{cor:inf}} that the semilocking boundary point is difficult to be checked as an $\omega$-$limit$ point, especially when it is a boundary equilibrium. In this subsection, we will provide some sufficient conditions to suggest persistence of MASs under the framwork of the Lyapunonv Function PDEs method. For this purpose, we begin with the following lemma.
	
\begin{lemma}\label{lema:4.1}
Given a MAS $(\mathcal{S,C,R},k)$ with the dynamics of {(\ref{eq:de1})}, let its Lyapunov Function PDEs (\ref{eq:pde}) and (\ref{eq:pdea}) admit a solution $f\in \mathscr{C}^{1}(\mathbb{R}^{\mathcal{S}}_{>0};\mathbb{R})$. Then for any semilocking boundary point $\bar{x}\in\mathbb{R}^{\mathcal{S}}_{\geq 0}$ we have 
    \begin{align*}
   & \varlimsup\limits_{x\rightarrow \bar{x}}\dot{f}(x)
    =\varlimsup\limits_{x\rightarrow \bar{x}}\sum_{i=1}^{r}k_{i}x^{v_{\cdot i}}(v'_{\cdot i}-v_{\cdot i})^{\top}\triangledown f(x)\\
& \leq\underbrace{\varlimsup\limits_{x\rightarrow \bar{x}}\sum_{\substack{\mathrm{supp}~v_{\cdot i}\subseteq \mathrm{supp} ~\bar{x}\\\mathrm{supp}~v'_{\cdot i}\subseteq \mathrm{supp} ~\bar{x}}}k_{i}x^{v_{\cdot i}}(v'_{\cdot i}-v_{\cdot i})^{\top}\triangledown f(x)}_{F_1(x)}+\underbrace{\varlimsup\limits_{x\rightarrow \bar{x}}\sum_{\substack{\mathrm{supp} ~v_{\cdot i}\subseteq \mathrm{supp} ~\bar{x}\\\mathrm{supp}~v'_{\cdot i}\nsubseteq \mathrm{supp} ~\bar{x}}}k_{i}x^{v_{\cdot i}}(v'_{\cdot i}-v_{\cdot i})^{\top}\triangledown f(x)}_{F_2(x)}\\
    &+\underbrace{\varlimsup\limits_{x\rightarrow \bar{x}}\sum_{\substack{\mathrm{supp}~v_{\cdot i}\nsubseteq \mathrm{supp}~\bar{x}\\\mathrm{supp}~v'_{\cdot i}\subseteq \mathrm{supp}~\bar{x}}}k_{i}x^{v_{\cdot i}}(v'_{\cdot i}-v_{\cdot i})^{\top}\triangledown f(x)}_{F_3(x)}    +\underbrace{\varlimsup\limits_{x\rightarrow \bar{x}}\sum_{\substack{\mathrm{supp}~v_{\cdot i}\nsubseteq \mathrm{supp}~\bar{x}\\\mathrm{supp}~v'_{\cdot i}\nsubseteq \mathrm{supp}~\bar{x}}}k_{i}x^{v_{\cdot i}}(v'_{\cdot i}-v_{\cdot i})^{\top}\triangledown f(x)}_{F_4(x)} \\
 & =F_1(x)+F_3(x)+F_4(x).   
    \end{align*}
    	\end{lemma}
	\noindent{\textbf{Proof.}}
	For the above equation, the first two pieces of equality hold obviously. For the third one, since $W\triangleq \mathrm{supp}^{c}\bar{x}$ is a semilocking set, there does not exist a reaction $v_{\cdot i}\to v'_{\cdot i}$ in $(\mathcal{S,C,R},k)$ such that $\mathrm{supp}~v_{\cdot i}\subseteq \mathrm{supp}~\bar{x}$ and $\mathrm{supp}~v'_{\cdot i}\nsubseteq \mathrm{supp}~\bar{x}$. Hence, $F_2(x)=0$, which completes the proof. $\Box$
	
The above result may simplify the computation of $\varlimsup\limits_{x\rightarrow \bar{x}}\dot{f}(x)$ when $\bar{x}$ is a semilocking boundary point. We then set some sufficient conditions to lead to the result of $\varlimsup\limits_{x\rightarrow \bar{x}}\dot{f}(x)$ falling within the range presented in \textbf{Theorem \ref{Th:4.5}}.  

\begin{theorem}\label{gaoThm4.6}
For a MAS $(\mathcal{S,C,R},k)$ governed by {(\ref{eq:de1})}, assume its Lyapunov Function PDEs (\ref{eq:pde}) and (\ref{eq:pdea}) admit a nonnegative solution $f\in \mathscr{C}^{2}(\mathbb{R}^{\mathcal{S}}_{>0};\mathbb{R})$. Then for any semilocking boundary non-equilibrium point $\bar{x}\in\mathbb{R}^{\mathcal{S}}_{\geq 0}$ of this network, $\bar{x}$ is not an $\omega$-$limit$ point if for any $x\in\mathbb{R}^{\mathcal{S}}_{> 0}$ the Hessian matrix of $f(x)$ is a positive definite diagonal matrix.
\end{theorem}
\noindent{\textbf{Proof.}}
See Appendix. $\Box$

Note that the results given in the above theorem are relatively weak, since only semilocking boundary non-equilibrium points can be captured. It is still unknown for semilocking boundary equilibrium points whether they are $\omega$-$limit$ points or not. A direct application of \textbf{Theorem \ref{gaoThm4.6}} is to complex balanced MASs.

\begin{corollary}\label{complexcase}
Assume a MAS $(\mathcal{S,C,R},k)$ described by {(\ref{eq:de1})} is complex balanced. Then any of its semilocking boundary non-equilibrium points $\bar{x}\in\mathbb{R}^{\mathcal{S}}_{\geq 0}$ is not an $\omega$-$limit$ point.
\end{corollary}
\noindent{\textbf{Proof.}}
It has been proved \cite{Fang} that the well-known pseudo-Helmholtz free energy function \cite{Horn1972}, given by 
\begin{equation}
	G(x)=\sum_{j=1}^n
	x_{j}\left(\ln x_{j}-\ln x^{*}_{j}-1\right)+x^{*}_{j},~~~x\in\mathbb{R}^{\mathcal{S}}_{>0},
	\end{equation}
 is a nonnegative solution of the corresponding Lyapunov Function PDEs, where $x^*\in\mathbb{R}^{\mathcal{S}}_{>0}$ is an equilibrium in the complex balanced system. Note that the eqilibrium $x^*$ should be chosen such that $x^*-\bar{x}\in\mathscr{S}$, which must be existing since each positive stoichiometric compatibility class must contain and only contains an equilibrium \cite{Horn1972}. Clearly, $G\in\mathscr{C}^2(\mathbb{R}^{\mathcal{S}}_{>0};\mathbb{R})$ and its Hessian matrix, written as
\begin{equation*}
	 \triangledown^{2}G(x)=
	 \left(\begin{matrix}
	 \frac{1}{x_{1}}& & \\
	  &\ddots& \\
	  & &\frac{1}{x_{n}}
	 \end{matrix}\right),
	\end{equation*}
is a positive definite diagonal matrix. The result is straightforward based on \textbf{Theorem \ref{gaoThm4.6}}. $\Box$

It should be pointed out that this result was also reached in \cite{Siegel} using the $\omega$-limit set theorem. We illustrate an application of \textbf{Theorem \ref{gaoThm4.6}} through the following example which is obviously not a complex balanced network even not a weakly reversible network.

\begin{example}
A MAS follows
   $$\xymatrix{ S_{1}\xrightarrow{k_{1}} S_{2},}
   \xymatrix{2S_{2}\xrightarrow{k_{2}}S_{2}+S_{1},}
  \xymatrix{2S_{1}+S_{3} \ar @{ -^{>}}^{~~~k_3}  @< 1pt> [r] &  2S_3\ar @{ -^{>}}^{~~~k_4} @< 1pt> [l] }.$$
By setting $k_{1}=k_{2}=k_{3}=k_{4}=1$, we have its dynamics to be
\begin{equation}\label{ex5dynamics}
\dot{x}(t)=\left(\begin{array}{rrrr}
        -1 & 1& -2& 2 \\
         1& -1& 0& 0\\
         0& 0& 1& -1\\
    \end{array}\right)\left(\begin{array}{c}
         x_{1} \\
         x_{2}^{2}\\
         x_{1}^{2}x_{3}\\
         x_{3}^{2}\\
    \end{array}\right)=\left(\begin{array}{c}
         -x_{1}+x_{2}^{2}-2x_{1}^{2}x_{3}+2x_3^2\\
        x_{1}-x_{2}^{2}\\
         x_{1}^{2}x_{3}-x_{3}^{2}
    \end{array}\right),
    \end{equation}
and the Lyapunov Function PDE {(\ref{eq:pde})} to be 
\begin{align*}
&x_1+x_{2}^{2}+x_{1}^{2}x_{3}+x_3^{2}-x_1\exp\left\{-\frac{\partial f}{\partial x_1}+\frac{\partial f}{\partial x_2}\right\}-x_{2}^{2}\exp\left\{\frac{\partial f}{\partial x_1}-\frac{\partial f}{\partial x_2}\right\}\\&-x_{1}^{2}x_{3}\exp\left\{-2\frac{\partial f}{\partial x_1}+\frac{\partial f}{\partial x_3}\right\}
-x_{3}^{2}\exp\left\{2\frac{\partial f}{\partial x_1}-\frac{\partial f}{\partial x_3}\right\}=0.
\end{align*}
We can verify that the function $$f(x)= x_1\ln{x_1}-x_1+2(x_{2}\ln{x_2}-x_2)+x_3\ln{x_3}-x_3+4$$ is a nonnegative solution of the above PDE. Its time derivative is
\begin{align*}
        \dot{f}(x)=(-x_1+x_{2}^{2}-2x_{1}^{2}x_{3}+2x_{3}^{2})\ln{x_{1}}+2(x_{1}-x_{2}^{2})\ln{x_{2}}+(x_{1}^{2}x_{3}-x_{3}^{2})\ln{x_{3}},
    \end{align*}
and Hessian matrix is $\triangledown^{2}f(x)=\mathrm{diag}(1/x_1,2/x_2,1/x_3)$ that is obviously a positive diagonal matrix in $\mathbb{R}^3_{>0}$. Hence, from \text{Theorem \ref{gaoThm4.6}} any semilocking boundary non-equilibrium point in this MAS is not an $\omega$-$limit$ point. 

We exhibit more details here. For simplicity, consider a representative semilocking boundary point $\bar{x}=(\bar{x}_{1},\bar{x}_2,0)^\top$, where $\bar{x}_1,\bar{x}_2>0$ and $W\triangleq \mathrm{supp}^{c}\bar{x}=\{S_{3}\}$, then we have
    \begin{align*}
           \lim\limits_{x\rightarrow\bar{x}} \dot{f}(x)&=\lim\limits_{x\rightarrow\bar{x}}\left[(-x_1+x_{2}^{2}-2x_{1}^{2}x_{3}+2x_{3}^{2})\ln{x_{1}}+2(x_{1}-x_{2}^{2})\ln{x_{2}}+(x_{1}^{2}x_{3}-x_{3}^{2})\ln{x_{3}}\right]\\&=-(\bar{x}_{1}-\bar{x}_{2}^{2})(\ln{\bar{x}_1}-2\ln{\bar{x}_{2}}).
    \end{align*}
 Therefore, if $\bar{x}_1\neq\bar{x}_2^2$, which suggests $\bar{x}$ to be not a boundary equilibrium from the dynamical equation {(\ref{ex5dynamics})}, then $\lim\limits_{x\rightarrow\bar{x}} \dot{f}(x)< 0$. $\bar{x}$ is not an $\omega$-$limit$ point. However, if $\bar{x}_1=\bar{x}_2^2$, then $\bar{x}$ is a boundary equilibrium at which $\lim\limits_{x\rightarrow\bar{x}} \dot{f}(x)= 0$. In this case, it is unknown whether $\bar{x}$ is an $\omega$-$limit$ point or not.    
\end{example}
\section{Persistence of $1d$ $\mathcal{W}_{I}$-endotactic networks}\label{sec:5}
In this section, we will use the proposed Lyapunov Function PDEs technique to prove that all of $1d~\mathcal{W}_{I}$-endotactic networks have persistence. 

\subsection{Solution of the Lyapunov Function PDEs of $1d$ MASs}
The key point to perform persistence analysis utilizing the Lyapunov Function PDEs technique lies in achieving the explicit expression of their solution. 

For a $1d$ MAS $(\mathcal{S,C,R},k)$, the Lyapunov Function PDEs can be simplified as
\begin{equation}\label{eq:1pde}
		(u-1)
		\left[\sum_{\{i|m_{i}>0\}}\left(k_{i}x^{v_{\cdot
				i}}\right)\left(\sum_{j=0}^{m_{i}-1}u^{j}\right)
		+\sum_{\{i|m_{i}<0\}}\left(k_{i}x^{v_{\cdot
				i}}\right)\left(-\sum_{j=m_{i}}^{-1}u^{j}\right)\right]=0
		\end{equation}
and {(\ref{eq:pdea})}, where $u=\exp\{b^{\top}\triangledown f(x)\},~b\in \mathbb{R}^{\mathcal{S}}\backslash \{\mathbbold{0}_{n}\}$ is a group of base of the stoichiometric subspace $\mathscr{S}$, and $m_{i}\in \mathbb{Z}\backslash\{0\},~i=1,~\cdots,~r$, satisfy
$v'_{\cdot i}-v_{\cdot i}=m_{i}b.$ For simplicity of symbols, denote a a scalar function $h(x,u)$ by
\begin{equation}\label{eq:gaoh}
	  h(x,u)=\sum_{\{i|m_{i}>0\}}(k_{i}x^{v_{\cdot i}})\left(\sum_{j=0}^{m_{i}-1}u^{j}\right)+\sum_{\{i|m_{i}<0\}}(k_{i}x^{v_{\cdot i}})\left(-\sum_{j=m_{i}}^{-1}u^{j}\right).
\end{equation}
It has been shown that there exists a unique $\tilde{u}\in \mathscr{C}^{2}(\mathbb{R}^{\mathcal{S}}_{>0};\mathbb{R})$ supporting $h(x,\tilde{u}(x))=0$ if there exists an positive equilibrium in this network\cite{Fang}.

\begin{lemma}[\cite{Fang}]\label{pro:4.5}
For a $1d$ MAS $(\mathcal{S,C,R},k)$ modeled by {(\ref{eq:de1})} admitting an positive equilibrium, the function defined by 
\begin{equation}\label{eq:1pdes}
   f(x)=\int_{0}^{\gamma(x)}\ln \tilde{u}(y^{\dag}(x)+\tau b)d\tau
\end{equation}
is a solution of its Lyapunov Function PDEs {(\ref{eq:1pde})} and {(\ref{eq:gaoh})}, where $y^{\dag}\in \mathscr{C}^{2}(\mathbb{R}^{\mathcal{S}}_{>0};\mathbb{R}^{\mathcal{S}}_{>0})$, $\gamma\in \mathscr{C}^{2}(\mathbb{R}^{\mathcal{S}}_{>0};\mathbb{R})$, constrained by
	\begin{equation*}
	x=y^{\dag}(x)+\gamma(x)b~~~ and ~~~\gamma(x+\delta b)=\gamma(x)+\delta,~~\forall \delta \in \mathbb{R},
	\end{equation*}
and $\tilde{u},b$ share the same meanings as in {(\ref{eq:1pde}) and (\ref{eq:gaoh})}.
	\end{lemma}

\textbf{Proposition \ref{pro:4.5}} makes it clear that there is an explicit solution for the Lyapunov function PDEs of any $1d$ network. Hence, it is possible to apply \textbf{Theorem \ref{Th:4.5}} to analyze persistence within this class of networks. Remarkably, the result given in {(\ref{Th:3.4})} implies that we only need to consider semilocking boundary points in the subsequent analysis.

\begin{proposition}\label{pro:5.2}
Each semilocking boundary point in a $1d$ MAS $(\mathcal{S,C,R},k)$ is a boundary equilibrium point.
\end{proposition}
\noindent{\textbf{Proof.}}
In a $1d$ MAS $(\mathcal{S,C,R},k)$, any reaction vector may be expressed as 
\begin{equation}\label{1drv}
 v'_{\cdot i}-v_{\cdot i}=m_{i}b,~~~\forall~ i=1,\cdots,r,
 \end{equation}
where $b\in\mathbb{R}^{n}$ $\backslash \{\mathbbold{0}_{n}\}$ is a group of base in $\mathscr{S}$ and $m_{i}\in \mathbb{Z}\backslash \{0\}$ is the corresponding coefficient. Let $\bar{x}\in\mathbb{R}^{\mathcal{S}}_{\geq 0}$ represent a semilocking boundary point in the $1d$ MAS under consideration, i.e., $I=\mathrm{supp}^{c}~\bar{x}$ is a semilocking set, then any reaction $v_{\cdot i}\to v'_{\cdot i}$ in $(\mathcal{S,C,R},k)$ does not support 
$\mathrm{supp}~v_{\cdot i}\subseteq \mathrm{supp}~\bar{x}$ and  $\mathrm{supp}~v'_{\cdot i}\nsubseteq \mathrm{supp}~\bar{x}$. 

Assume there includes a reaction $v_{\cdot l}\to v'_{\cdot l}$ in $(\mathcal{S,C,R},k)$ supporting $\text{supp}~v_{\cdot l}\subseteq \text{supp}~ \bar{x}$ and $\text{supp}~v'_{\cdot l}\subseteq \text{supp} ~\bar{x}$. Then, $\forall ~j\in I$ we have $v'_{jl}-v_{jl}=0$, which together with {(\ref{1drv})} suggests $b_j=0$. This result conversely renders that there is no change on the concentration of species $S_j$, $\forall ~j\in I$, within all reactions in $(\mathcal{S,C,R},k)$. We may naturally ignore these species in persistence analysis since they are requested to have nonzero concentrations in the beginning and will not affect the persistence property of system. Hence, we only need to be concerned with those species whose concentrations will change within reactions. We suppose all species concentrations will vary during the reaction process, i.e., every entry of $b$ is not zero. This means that among all reactions it is impossible to have a reaction  $v_{\cdot l}\to v'_{\cdot l}$ supporting $\text{supp}~v_{\cdot l}\subseteq \text{supp}~ \bar{x}$ and $\text{supp}~v'_{\cdot l}\subseteq \text{supp} ~\bar{x}$.  

Based on the above results obtained, every reaction $v_{\cdot i}\to v'_{\cdot i}$ in a $1d$ MAS does not support $\mathrm{supp}~v_{\cdot i}\subseteq \mathrm{supp}~\bar{x}$. The network can only include reactions satisfying $\mathrm{supp}~v_{\cdot i} \nsubseteq \mathrm{supp}~\bar{x}$. As a result, $I$ is a locking set as well as a semilocking set. Therefore, $\bar{x}$ is a boundary equilibrium point. $\Box$

\begin{remark}\label{remk:gao}
The proof process of \textbf{Proposition \ref{pro:5.2}} indicates that we only need to select a group of base $b$ in $\mathscr{S}$ with every entry not to be zero when we make persistence analysis for $1d$ MAS. We thus limit $b\in\mathbb{R}^{\mathcal{S}}_{\neq 0}$ in what follows.
\end{remark}

$\textbf{Proposition \ref{pro:5.2}}$ tells us that in the case of $1d$ MASs, the concerning semilocking boundary points are all boundary equilibria, at which the solution of the corresponding Lyapunov function PDEs will stop change. This will invalidate the first two conditions given in \textbf{Theorem \ref{Th:4.5}} for suggesting persistence. We thus switch our attention to the application of the third condition, i.e., local maximum condition at the boundary point, and consider $1d~\mathcal{W}_{I}$-endotactic networks, specially.

\subsection{Persistence analysis of $1d~\mathcal{W}_{I}$-endotactic networks}
We firstly give the definition of $\mathcal{W}_{I}$-endotactic networks. 
\begin{definition}\label{WInetwork}
Given a reaction network $(\mathcal{S,C,R})$, let $I$ represent a semilocking set attached to this network. If $w_{I}\in\{0,1\}^{\mathcal{S}}$ is defined by
$$ w_{I_j}=
\begin{cases}
       &1,~~ j\in I,  \\
     &0, ~~j\notin I,
\end{cases}
 $$
and $\mathcal{W}_{I}=\{w_I|I~ is~ any~ semilocking ~set~ attached~ to~ (\mathcal{S,C,R}).\}$, then the network $(\mathcal{S,C,R})$ is called $\mathcal{W}_{I}$-endotactic if for every $w_I\in\mathcal{W}_I$ this network is $w_{I}$-endotactic.  
\end{definition}

Comparisons of this definition to \textbf{Definition \ref{def:Endotatic}} might reveal that an endotatic network must be $\mathcal{W}_{I}$-endotactic, and a $\mathcal{W}_{I}$-endotactic network must be $w_I$-endotactic. However, the inverse propositions are not ture. An example of $\mathcal{W}_{I}$-endotactic networks is given as follows.

\begin{figure}
        \centering
        \begin{minipage}[c]{0.45\textwidth}
        \centering
          \includegraphics[width=1\textwidth]{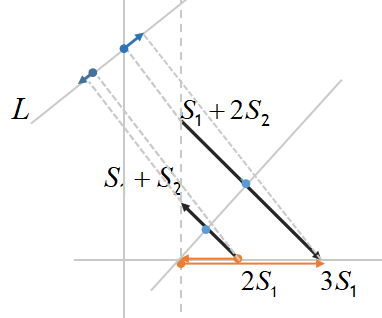}
        \end{minipage}
         \caption{A schematic diagram of network in \textbf{Example \ref{ex:z2}}, and the corresponding projections on some straight lines.}
        \label{fig:3}
    \end{figure}

\begin{example}\label{ex:z2}
For a network like 
\begin{align*}
      2S_{1}&\to S_{1}+S_{2},\\
      2S_{2}+S_{1}&\to 3S_{1},
\end{align*}
it can be seen from \textbf{Fig. \ref{fig:3}} that this network is not endotactic since its projection onto the straight line $L$ fails to support the endotactic property. However, we can easily verify it to be $\mathcal{W}_{I}$-endotactic. The entire semilocking sets attached to this network are $\{S_1\}$ and $\{S_1,S_2\}$, then we get  $w_{\{S_1\}}=(1,0)$ and $w_{\{S_1,S_2\}}=(1,1)$. Also from \textbf{Fig. \ref{fig:3}}, the projected networks onto the lines where $w_{\{S_1\}}$ and $w_{\{S_1,S_2\}}$ lie, respectively, both show to be $w_I$-endotactic. So this network is $\mathcal{W}_{I}$-endotactic.
\end{example}

A $\mathcal{W}_{I}$-endotactic network $(\mathcal{S,C,R})$ is a $1d$ $\mathcal{W}_{I}$-endotactic network if $\mathrm{dim}\mathscr{S}=1$, which is our main concern later. Note that the persistence criteria by means of the $\omega$-$limit$ point analysis, such as those presented in \textbf{Theorem \ref{Th:4.5}}, request that the trajectories in the MAS are bounded. The following property for $1d$ MASs will guarantee the boundedness of trajectories in $1d~w_I$-endotactic networks.  
  
\begin{lemma}\label{lm:09}
Consider a $1d$ MAS $(\mathcal{S,C,R},k)$. Let $x(0)\in\mathbb{R}^{\mathcal{S}}_{\geq 0}$ be any initial point of the system. Then for arbitrary trajectory $x(t)\in\mathbb{R}^{\mathcal{S}}_{\geq 0}$ originating from $x(0)$, if the concentration of one species tends to infinity as time runs to infinity, so do all species concentrations. 
\end{lemma}
\noindent{\textbf{Proof.}}
From $\mathrm{dim}\mathscr{S}=1$, let $b\in\mathbb{R}^{\mathcal{S}}_{\neq 0}$, as suggested in \textbf{Remark \ref{remk:gao}}, be a group of base in this $1d$ MAS, then there exists a real function $a(t)$ satisfying $x(t)-x(0)=a(t)b$. Also, denote the positive/negative terms index set of $b$ by $P_b$/$N_b$. Without loss of generality, suppose $l\in P_b$ and $\lim\limits_{t\rightarrow \infty}x_{l}(t)=\infty$, which means $\lim\limits_{t\rightarrow \infty}a(t)=\infty$, so $\forall~j\in P_b$ we have $\lim\limits_{t\rightarrow \infty}x_{j}(t)=\infty$ while $\forall~j\in N_b$ we have $\lim\limits_{t\rightarrow \infty}x_{j}(t)=-\infty$. The latter is apparently not true. At this juncture, we obtain $N_b=\varnothing$. Therefore, either $P_b$ or $N_b$ must be empty under the given conditions. The result is thus true.$\Box$

Then we can catch the boundedness of trajectories in $1d~w_{\mathcal{S}}$-endotactic networks.

\begin{lemma}\label{lm:bd}
Any trajectory in a $1d~w_{\mathcal{S}}$-endotactic MAS $(\mathcal{S,C,R},k)$ is bounded.
\end{lemma}
\noindent{\textbf{Proof.}}
We use the proof by contradiction to reach the result. Assume a trajectory $x(t)$ in the $1d~w_{\mathcal{S}}$-endotactic MAS $(\mathcal{S,C,R},k)$ is unbounded as time runs to infinity. From \textbf{Lemma \ref{lm:09}}, we have known every species concentration $x_{j}(t)$ subject to $\lim\limits_{t\rightarrow \infty}x_{j}(t)=\infty,~\forall~S_j\in\mathcal{S}$. Using the same variables as used in the proof of \textbf{Lemma \ref{lm:09}}, like $x(0)$, $b$, $a(t)$, etc., and setting all entries in $b$ to be positive for simplicity but without loss of generality, we get  
 	  \begin{equation*}
	     \lim\limits_{t\rightarrow \infty}\frac{x_{j}(t)}{x_{1}(t)}= \lim\limits_{t\rightarrow \infty}\frac{x_{j}(0)+a(t)b_{j}}{x_{1}(0)+a(t)b_{1}}=\frac{b_{j}}{b_{1}}>0, ~\forall~ S_j\in\mathcal{S}.
	    \end{equation*}
This means $\exists~ T_1,\beta_j,\bar{\beta}_j >0$ such that when $t>T_1$, there is 
\begin{equation*}
	      \beta_{j}\leq\frac{x_{j}(t)}{x_{1}(t)}\leq \bar{\beta}_{j},~\forall~ S_j\in\mathcal{S}. 
	  \end{equation*}
i.e., 	  
\begin{equation}\label{eq:wmax}
	      \beta^{v_{\cdot i}}x_1^{\langle v_{\cdot i}, w_{\mathcal{S}} \rangle}\leq x^{v_{\cdot i}}\leq {\bar{\beta}}^{v_{\cdot i}}x_1^{\langle v_{\cdot i}, w_\mathcal{S} \rangle},~\forall~v_{\cdot i}\in\mathcal{C}
	  \end{equation}
where $\beta=(\beta_1,\cdots,\beta_n)^\top$, $\bar{\beta}=(\bar{\beta}_1,\cdots,\bar{\beta}_n)^\top$, $\mathcal{S}$ is naturally a semilocking set, and  $w_\mathcal{S}=(1,\cdots,1)^{\top}$ following \textbf {Definition \ref{WInetwork}}. 
Further define the set composed of all $\leq_{w_{\mathcal{S}}}$-maximal reactant complexes by $\mathcal{Y}$, i.e., $\forall~ v_{\cdot l}\in \mathcal{Y}~\mathrm{and}~\forall~ v_{\cdot i}\in \mathcal{C}\backslash\mathcal{Y}$, there is $\langle v_{\cdot l}, w_{\mathcal{S}}\rangle > \langle v_{\cdot i}, w_{\mathcal{S}}\rangle$.
Then utilizing the explosiveness of the exponential function, and the fact that $\beta^{v_{\cdot l}}$ and ${\bar{\beta}}^{v_{\cdot i}}$ are constants, we get $\exists~T_2>0$ such that when $t>T_2$ there is $${\bar{\beta}}^{v_{\cdot i}}x_1^{\langle v_{\cdot i}, w_\mathcal{S} \rangle}<\beta^{v_{\cdot l}}x_1^{\langle v_{\cdot l}, w_\mathcal{S} \rangle}.$$
Combining the above inequality with {(\ref{eq:wmax})} yields that when $t>\max{\{T_1,T_2\}}$, we have
$$x^{v_{\cdot i}}\leq {\bar{\beta}}^{v_{\cdot i}}x_1^{\langle v_{\cdot i}, w_\mathcal{S} \rangle}<\beta^{v_{\cdot l}}x_1^{\langle v_{\cdot l}, w_\mathcal{S} \rangle}\leq x^{v_{\cdot l}}, ~~~\forall~ v_{\cdot l}\in \mathcal{Y}, v_{\cdot i}\in \mathcal{C}\backslash\mathcal{Y}.$$
For the same reasons as above, there must exist $T_3>0$ such that when $t>\max{\{T_1,T_2,T_3\}}$ we have
\begin{equation*}
k_{l}x^{v_{\cdot l}}>k_{i}x^{v_{\cdot i}}, \forall~ v_{\cdot l}\in \mathcal{Y}, v_{\cdot i}\in \mathcal{C}\backslash\mathcal{Y},
\end{equation*}
and $\exists~T_4>0$ such that when $t>\max{\{T_1,T_2,T_3,T_4\}}$ we further have
\begin{equation}\label{eq:11}
\sum_{v_{\cdot i}\in \mathcal{Y}}k_{i}x^{v_{\cdot i}}>\sum_{v_{\cdot i}\in\mathcal{C}\backslash\mathcal{Y}}k_{i}x^{v_{\cdot i}}.
\end{equation}
Note that the network under consideration is $1d~w_{\mathcal{S}}$-endotactic, and also $w_{\mathcal{S}}$ is not orthogonal to $\mathscr{S}$ (due to $\langle w_{\mathcal{S}}, b\rangle \neq 0$), from \textbf{Definition \ref{def:Endotatic}} we thus get 
$$\langle v'_{\cdot i}-v_{\cdot i}, w_{\mathcal{S}}\rangle<0,\forall~ v_{\cdot i}\in \mathcal{Y},$$ which combines {(\ref{1drv})} to yield 
$m_i<0$. 

The above results are used to rewrite the dynamics of the network, which gives
\begin{align*}
\dot{x}(t)&=\sum_{i=1}^{r}k_{i}x^{v_{\cdot i}}(v'_{\cdot i}-v_{\cdot i}) \\
&=\sum_{\{i|v_{\cdot i}\notin\mathcal{Y}\}}k_{i}x^{v_{\cdot i}}m_{i}b+\sum_{\{i|v_{\cdot i}\in\mathcal{Y}\}}k_{i}x^{v_{\cdot i}}m_{i}b\\
&=\sum_{\{j|m_{i}>0\}}k_{i}x^{v_{\cdot i}}m_{i}b-\sum_{\{i|m_{i}<0\}}k_{i}x^{v_{\cdot i}}|m_{i}|b.
\end{align*}
Again, by using the explosiveness of the exponential function and {(\ref{eq:11})}, we obtain that $\exists~T_4>0$ satisfying when $t>\max{\{T_1,T_2,T_3,T_4\}}$, every entry of $\sum_{\{j|m_{i}<0\}}k_{i}x^{v_{\cdot i}}|m_{i}|b$ is larger than that of $\sum_{\{i|m_{i}>0\}}k_{i}x^{v_{\cdot i}}m_{i}b$. As a result, there is $\dot{x}_j(t)<0$ if $t>\max{\{T_1,T_2,T_3,T_4\}}$, $\forall~S_j$. This is contrary to the assumption that $x(t)$ is unbounded as $t\to \infty$. Therefore, the result holds. 

Note that \textbf{Proposition \ref{pro:4.5}} holds on the premise of the existence of a positive equilibrium in the $1d$ MAS. The role is to ensure the terms in {(\ref{eq:gaoh})} neither $\{i|m_{i}>0\}$ nor $\{i|m_{i}<0\}$ is empty \cite{Fang}, which will further guarantee that there exists a unique $\tilde{u}\in \mathscr{C}^{2}$ supporting $h(x,\tilde{u}(x))=0$. Clearly, $1d~\mathcal{W}_I$ endotactic networks admit the result that neither $\{i|m_{i}>0\}$ nor $\{i|m_{i}<0\}$ is empty. Hence, for a $1d~\mathcal{W}_I$ endotactic network, $\tilde{u}\in \mathscr{C}^{2}$ is unique to support $h(x,\tilde{u}(x))=0$, and the function $f(x)$ given in {(\ref{eq:1pde})} is the solution of its Lyapunov function PDEs, no matter whether the system admits a positive equilibrium. We further analyze the property of the solution of the Lyapunov function PDEs for a $1d~\mathcal{W}_{I}$-endotactic MAS, which will approach our main result.
\begin{lemma}\label{lm:bd2}
For a $1d~\mathcal{W}_{I}$-endotactic MAS $(\mathcal{S,C,R},k)$, let $\bar{x}\in\mathbb{R}^{\mathcal{S}}_{\geq 0}$ represent any semilocking boundary point, and $x_0\in\mathscr{S}(\bar{x})\cap \mathbb{R}^{\mathcal{S}}_{> 0}$ represent any initial point in this MAS. Without loss of generality, assume that $I\triangleq \mathrm{supp}~\bar{x}$ belongs to the positive terms index set of the base $b$ in $\mathscr{S}$. Then the solution of its Lyapuonv function PDEs {(\ref{eq:1pde})}, given by {(\ref{eq:1pdes})}, is constrained by $\lim\limits_{\delta\rightarrow 0^{+}}\tilde{u}(\bar{x}+\delta b)<1.$  
\end{lemma}
\noindent{\textbf{Proof.}}
Based on \textbf{Proposition \ref{pro:4.5}}, the positive function $\tilde{u}(x)$ emerging in the solution $f(x)$ of {(\ref{eq:1pdes})} should satisfy $h(x,\tilde{u}(x))=0$, i.e., {(\ref{eq:gaoh})} becomes
\begin{equation}\label{compaGao}
\sum_{\{i|m_{i}>0\}}(k_{i}x^{v_{\cdot i}})\left(\sum_{j=0}^{m_{i}-1}\tilde{u}^{j}(x)\right)=\sum_{\{i|m_{i}<0\}}(k_{i}x^{v_{\cdot i}})\left(\sum_{j=m_{i}}^{-1}\tilde{u}^{j}(x)\right).
\end{equation}
It is clear that $I\triangleq \mathrm{supp}^c~\bar{x}$ is a semilocking set. Without loss of generality, $I$ is set as $\{S_{1},\cdots,S_{n_{1}}\}$ in the following proof, where $1\leq n_{1}\leq n$, then $\bar{x}=(\underbrace{0,\cdots,0}_{n_{1}},\bar{x}_{n_1+1},\cdots,\bar{x}_n)^{\top}$ and 
$w_{I}=(\underbrace{1,\cdots,1}_{n_{1}},0,\cdots,0)^{\top}$ according to \textbf{Definition \ref{WInetwork}}. Also, we define a set $\mathcal{Z}$ of $\leq_{w_{I}}$-minimal reactant complexes, i.e., $\langle v_{\cdot l}-v_{\cdot i}, w_I \rangle \leq 0,~\forall~ v_{\cdot i}\in \mathcal{C}~\mathrm{and}~\forall~v_{\cdot l}\in\mathcal{Z}.$

Since $I\triangleq \mathrm{supp}^c~\bar{x}$ belongs to the positive terms index set of the base $b$, $\exists~\bar{a}>0$ such that $\bar{x}=x(0)-\bar{a}b$. For those $x(t)\in\mathscr{S}(\bar{x} )$ near $\bar{x}$, there exists a real function $a(t)$ such that $x(t)=x(0)-a(t)b$. We thus get
\begin{align*}
  \lim_{a(t)\rightarrow \bar{a}^{-}}\frac{x_{j}(t)}{x_{1}(t)}=\frac{b_{j}}{b_{1}}, ~\forall~S_j\in I \mathrm{~and~}\lim_{a(t)\rightarrow \bar{a}^{-}}x_{j}=\bar{x}_{j},  ~\forall~S_j\in\mathcal{S}\backslash I.
\end{align*}
The first equation means that $\exists~\tau_1,\theta_j,\bar{\theta}_j>0$ such that when $\bar{a}-a(t)<\tau_{1}$ there exists
\begin{equation*}
    \theta_{j}\leq\frac{x_{j}(t)}{x_{1}(t)}\leq \bar{\theta}_{j}, ~~~\forall~S_j\in I.
\end{equation*} 
Therefore, for any reaction $v_{\cdot }\to v'_{\cdot i}$ in $(\mathcal{S,C,R},k)$ we have
\begin{align*}
\theta x_{1}^{\langle v_{\cdot i}, w_{I}\rangle}\leq\prod_{S_j\in I}x_{j}^{v_{ji}}\leq \bar{\theta}x_{1}^{\langle v_{\cdot i}, w_{I}\rangle},
\end{align*}
where $\theta=\prod_{S_j\in I}\theta_{j}^{v_{ji}}$ and $\bar{\theta}=\prod_{S_j\in I}\bar{\theta}_{j}^{v_{ji}}$ are two positive constants. Utilizing the property of the exponential function and $\mathcal{Z}$ being the set of $\leq_{w_{I}}$-minimal reactant complexes, we know that $\forall~v_{\cdot l}\in\mathcal{Z},~v_{\cdot i}\notin\mathcal{Z}$, $\exists~\tau_2>0$ such that when $\bar{a}-a(t)<\mathrm{min}\{\tau_{1},\tau_{2}\}$, i.e., $x_{1}<\mathrm{min}\{x_{1}(0)-(\bar{a}-\tau_{1})b_{1},x_{1}(0)-(\bar{a}-\tau_{2})b_{1}\}$, there is
$$\prod_{S_j\in I}x_{j}^{v_{jl}}\geq \theta x_{1}^{\langle v_{\cdot l}, w_{I}\rangle}>\bar{\theta} x_{1}^{\langle v_{\cdot i}, w_{I}\rangle}\geq \prod_{S_j\in I}x_{j}^{v_{ji}}.$$
This together with the facts that $x_j$ ($S_j\notin I$) is positively bounded and every reaction rate coefficient is a positive constant, means that $\exists~\tau_3>0$ such that when $\bar{a}-a(t)<\mathrm{min}\{\tau_{1},\tau_{2},\tau_3\}$, the inequality $k_{l}x^{v_{\cdot l}}>k_{i}x^{v_{\cdot i}}$ holds. For the same reason, $\exists~\tau_4>0$ supporting
\begin{equation}\label{eq:133}
\sum_{\{i|v_{\cdot i}\in \mathcal{Z}\}}k_{i}x^{v_{\cdot i}}>\sum_{\{i|v_{\cdot i}\notin \mathcal{Z}\}}k_{i}x^{v_{\cdot i}},
\end{equation}
when $\bar{a}-a(t)<\mathrm{min}\{\tau_{1},\tau_{2},\tau_3,\tau_4\}$.

Again, since $I\triangleq \mathrm{supp}^c~\bar{x}$ belongs to the positive terms index set of the base $b$, $w_{I}$ is not orthogonal to $\mathscr{S}$, i.e., $w_{I}$ not orthogonal to any reaction vector $v'_{\cdot i}-v_{\cdot i}$. We then get 
$$\langle v'_{\cdot i}-v_{\cdot i}, w_{I}\rangle >0,~~\forall~ v_{\cdot i}\in \mathcal{Z},$$
which is further combined with $v'_{\cdot i}-v_{\cdot i}=m_{i}b$ to yield $m_i>0,~\forall~v_{\cdot i}\in\mathcal{Z}.$ Therefore, {(\ref{eq:133})} may be rewritten as
\begin{align}\label{compaGao1}
\sum_{\{i|m_i>0\}}(k_{i}x^{v_{\cdot i}})>\sum_{\{i|m_i<0\}}(k_{i}x^{v_{\cdot i}}).
\end{align}
Therefore, in order to make {(\ref{compaGao})} hold, there must exist $\tilde{u}(x)<1$ when $x\in\mathbb{R}^{\mathcal{S}}_{>0}$ and $\bar{a}-a(t)$ is enough small.
This implies $\lim\limits_{\delta\rightarrow 0^{+}}\tilde{u}(\bar{x}+\delta b)<1$.

Finally, we present our main result as follows. 
\begin{theorem}\label{finalgao}
Each $1d~\mathcal{W}_{I}$-endotactic MAS $(\mathcal{S,C,R},k)$ is persistent.
\end{theorem}
\noindent{\textbf{Proof.}}
Let $\bar{x}\in\mathbb{R}^{\mathcal{S}}_{\geq 0}$ be any semilocking boundary point in the MAS under consideration (Only this case needs to be considered, since the case of non-semilocking boundary points has been addressed in \textbf{Theorem \ref{Th:3.4}}. Define $I\triangleq \mathrm{supp}^c~\bar{x}$ and $w_I$ following \textbf{Definition \ref{WInetwork}}, so $I$ is a semilocking set and the network in study is a $1d$ $w_I$-endotactic network, in which the trajectory originating any initial point is naturally bounded, as suggested in \textbf{Lemma \ref{lm:bd}}. We thus can apply \textbf{Theorem \ref{Th:4.5}} to analyze persistence for this MAS. The function $f(x)$ given in {(\ref{eq:1pdes})} is a solution of its Lyapunov function PDEs.

For any $x\in\mathscr{S}(\bar{x})\cap\mathbb{R}^{\mathcal{S}}_{>0}$ near $\bar{x}$, there must exist $\delta>0$ such that $x=\bar{x}+\delta b$, where $b$ is the base of $\mathscr{S}$, and is set with its positive terms index set including $I$, the same as what has been done in \textbf{Lemma \ref{lm:bd2}}. For convenient usage, we denote the functions in {(\ref{eq:1pdes})} evaluated at $x\to \bar{x}$ by
$$\tilde{u}(\bar{x})=\lim\limits_{ \delta\rightarrow 0^{+}}\tilde{u}(\bar{x}+\delta b) \mathrm{~and~} f(\bar{x})=\lim\limits_{\delta\rightarrow 0^{+}}f(\bar{x}+\delta b).$$ Then we have
\begin{align*}
\lim\limits_{\delta\rightarrow 0^{+}}\frac{f(\bar{x}+\delta b)-f(\bar{x})}{\delta}&=\lim\limits_{\delta\rightarrow 0^{+}}\frac{1}{\delta}\int_{\gamma(\bar{x})}^{\gamma(\bar{x}+\delta b)}\ln\tilde{u}(y^{\dag}(\bar{x})+\tau b)d\tau\\
&=\lim\limits_{\delta\rightarrow 0^{+}}\frac{1}{\delta}\int_{\gamma(\bar{x})}^{\gamma(\bar{x})+\delta}\ln\tilde{u}(y^{\dag}(\bar{x})+\tau b)d\tau\\
&=\ln \tilde{u}(y^{\dag}(\bar{x})+\gamma(\bar{x})b)\\
&=\ln\tilde{u}(\bar{x}).
\end{align*} 
From \textbf{Lemma \ref{lm:bd2}} that says $0<\tilde{u}(\bar{x})<1$, we get $$\lim\limits_{\delta\rightarrow 0^{+}}\frac{f(\bar{x}+\delta b)-f(\bar{x})}{\delta}=\lim\limits_{\delta\rightarrow 0^{+}}\frac{f(x)-f(\bar{x})}{\delta}<0.$$ Hence, $\bar{x}$ is a local maximum point for the function $f(x)$. Further based on \textbf{Theorem \ref{Th:4.5}}, we obtain the result that each $1d~\mathcal{W}_{I}$-endotactic MAS $(\mathcal{S,C,R},k)$ is persistent. $\Box$

\begin{example}\label{ex:z3}
We use \textbf{Example \ref{ex:z2}} to exhibit an application of \textbf{Theorem \ref{finalgao}}. As said in \textbf{Example \ref{ex:z2}}, the network is $\mathcal{W}_I$-endotactic. Also, from its stoichiometric matrix 
\begin{equation*}
    \varGamma=\left(\begin{array}{rr}
       -1  & 2 \\
        1 & -2
    \end{array}\right),
\end{equation*}
the network is a $1d~\mathcal{W}_I$-endotactic network. We assign $k_1,k_2>0$ to the first and second reaction in this network as the reaction rate coefficient, respectively. For simplicity, select $b=(1,-1)^{\top}$ as the base in $\mathscr{S}$, thus we get $m_{1}=-1,~ m_{2}=2$. 

There are three kinds of boundary points in this MAS, which are discussed in separate case.

$\mathrm{(1)}$ $\bar{x}=(\bar{x}_{1},0)^{\top}$, where $\bar{x}_1>0$. The set $W\triangleq \mathrm{supp}^c~\bar{x}=\{S_{2}\}$, which is obviously a non-semilocking set. According to \textbf{Theorem \ref{Th:3.4}}, this kinds of boundary points are not $\omega$-limit points.

$\mathrm{(2)}$ $\bar{x}=(0,\bar{x}_{2})^{\top}$, where $\bar{x}_2>0$. The set $I\triangleq \mathrm{supp}^c~\bar{x}=\{S_{1}\}$ is a semilocking set, and $b_1=1>0$, $w_{S_1}=(1,0)^{\top}$. The scalar function $h(x,u)$ in {(\ref{eq:gaoh})} follows 
$$h(x,u(x))=k_{2}x_{1}x_{2}^{2}(1+u(x))-k_{1}x_{1}^{2}u^{-1}(x).$$
We thus get $$\tilde{u}(x)=\frac{-k_{2}x_{2}+\sqrt{k_2^2x_{2}^{2}+4k_{1}k_2x_{1}}}{2k_{2}x_{2}}.$$
Further from 
$$\lim\limits_{\delta\rightarrow 0^{+}}\frac{f(\bar{x}+\delta b)-f(\bar{x})}{\delta}=\ln{\tilde{u}(\bar{x})}<0,$$
we conclude $\bar{x}=(0,\bar{x}_{2})^{\top}$ is not an $\omega$-limit point.

$\mathrm{(3)}$ $\bar{x}=(0,0)^{\top}$. The set $I\triangleq \mathrm{supp}^c~\bar{x}=\{S_{1},S_{2}\}$ is a semilocking set, and $w_{\{S_1,S_2\}}=(1,1)^\top$ is orthogonal to the stoichiometric subspace. It is impossible for this boundary point to belong to $\mathscr{S}(x(0))$ if $x(0)\in\mathbb{R}^{\mathcal{S}}_{>0}$. Namely, the origin is unattainable.

In short, this $1d$ $\mathcal{W}_{I}$-endotactic network is persistent, and moreover, the persistence is independent of the reaction rate coefficients. \end{example}
\section{Conclusions}\label{sec:6}
    This paper contributes to analyzing persistence of CRNs through the Lyapunov Function PDEs technique. Some new checkable criteria for persistence have been proposed by means of the solution, in the form of either its time derivative or itself, of the PDEs generated by the corresponding MAS. These criteria are successfully applied to capturing persistence of networks admitting only non-semilocking boundary points. However, for semilocking boundary points, only when they are also boundary non-equilibrium points and some moderate conditions are added, the proposed method can be valid. We further consider $1d$ MASs for which the Lyapunov Function PDEs are solvable \cite{Fang}, and define $\mathcal{W}_I$-endotactic networks, a set larger than endotactic networks but smaller than $w$-endotactic networks. We have proved that any $1d~\mathcal{W}_I$-endotactic network is persistent. Although this result is not completely new compared with those emerging in literature, such as the persistence of $1d$ strongly endotactic network covered by \cite{G2014}, the proposed Lyapunov Function PDEs technique provides an alternative to capture persistence of CRNs, and is exhibiting large potential. 
\begin{appendix}
\section{Proof of {Theorem \ref{gaoThm4.6}}}
Since $\bar{x}\in\mathbb{R}^{\mathcal{S}}_{\geq 0}$ is a semilocking boundary non-equilibrium point in the MAS, from \textbf{Lemma \ref{lema:4.1}} we have
	\begin{equation*}
		\varlimsup_{x\rightarrow \bar{x}}\dot{f}(x)\leq F_1(x)+F_3(x)+F_4(x).
		\end{equation*}
In addition, we have
	   \begin{align}
	   \lim_{x\rightarrow \bar{x}}\sum_{\mathrm{supp}~v_{\cdot i}\subseteq \mathrm{supp}~\bar{x}}k_{i}x^{v_{\cdot i}}&=\lim_{x\rightarrow \bar{x}}\sum_{\substack{\mathrm{supp}~ v_{\cdot i}\subseteq \mathrm{supp} ~\bar{x}\\\mathrm{supp}~v'_{\cdot i}\subseteq \mathrm{supp} ~\bar{x}}}k_{i}x^{v_{\cdot i}}+\lim_{x\rightarrow \bar{x}}\sum_{\substack{\mathrm{supp}~ v_{\cdot i}\subseteq \mathrm{supp}~ \bar{x}\\\mathrm{supp} ~v'_{\cdot i}\nsubseteq \mathrm{supp}~ \bar{x}}}k_{i}x^{v_{\cdot i}}\\
	   &=\lim_{x\rightarrow \bar{x}}\sum_{\substack{\mathrm{supp}~ v_{\cdot i}\subseteq \mathrm{supp} ~\bar{x}\\\mathrm{supp}~v'_{\cdot i}\subseteq \mathrm{supp}~ \bar{x}}}k_{i}x^{v_{\cdot i}}\\
	   &=\lim_{x\rightarrow \bar{x}}\sum_{\substack{\mathrm{supp}~ v_{\cdot i}\subseteq \mathrm{supp} ~\bar{x}\\\mathrm{supp}~v'_{\cdot i}\subseteq \mathrm{supp} ~\bar{x}}}k_{i}x^{v_{\cdot i}}+\lim_{x\rightarrow \bar{x}}\sum_{\substack{\mathrm{supp}~ v_{\cdot i}\nsubseteq \mathrm{supp}~ \bar{x}\\\mathrm{supp} ~v'_{\cdot i}\subseteq \mathrm{supp}~ \bar{x}}}k_{i}x^{v_{\cdot i}}\\
	   &=\lim_{x\rightarrow \bar{x}}\sum_{\mathrm{supp}~ v'_{\cdot i}\subseteq \mathrm{supp}~ \bar{x}}k_{i}x^{v_{\cdot i}},
	   \end{align}
which together with the boundary condition (a) {(\ref{eq:pdea})} yields
	   \begin{equation*}
	   \lim_{x\rightarrow \bar{x}}\sum_{\mathrm{supp}~ v'_{\cdot i}\subseteq \mathrm{supp}~ \bar{x}}k_{i}x^{v_{\cdot i}}=\lim_{x\rightarrow \bar{x}}\sum_{\mathrm{supp}~ v'_{\cdot i}\subseteq \mathrm{supp}~ \bar{x}}k_{i}x^{v_{\cdot i}}\exp\{(v'_{\cdot i}-v_{\cdot i})^{\top}\triangledown f(x)\}
	   \end{equation*}
By performing the Taylor expansion of $\exp\{(v'_{\cdot i}-v_{\cdot i})^{\top}\triangledown f(x)\}$ with respect to zero, we get
	   \begin{equation*}
	   \lim_{x\rightarrow \bar{x}}\left(\sum_{\mathrm{supp}~ v'_{\cdot i}\subseteq \mathrm{supp}~ \bar{x}}k_{i}x^{v_{\cdot i}}(v'_{\cdot i}-v_{\cdot i})^{\top}\triangledown f(x)+\sum_{\mathrm{supp}~ v'_{\cdot i}\subseteq \mathrm{supp}~ \bar{x}}k_{i}x^{v_{\cdot i}}\frac{e^{\alpha_{i}}}{2}[(v'_{\cdot i}-v_{\cdot i})^{\top}\triangledown f(x)]^2\right)=0,
	   \end{equation*}
	   where $\alpha_{i}\in \mathbb{R}$ lies between 0 and $(v'_{\cdot i}-v_{\cdot i})^{\top}\triangledown f(x).$ 
Further, we obtain
	     \begin{align*}
	        \varlimsup_{x\rightarrow \bar{x}}\sum_{\mathrm{supp}~ v'_{\cdot i}\subseteq \mathrm{supp}~ \bar{x}}k_{i}x^{v_{\cdot i}}(v'_{\cdot i}-v_{\cdot i})^{\top}\triangledown f(x)&\leq\underbrace{-\varlimsup\limits_{x\rightarrow\bar{x}}\sum_{\substack{\mathrm{supp}~ v_{\cdot i}\subseteq \mathrm{supp} ~\bar{x}\\\mathrm{supp}~v'_{\cdot i}\subseteq \mathrm{supp}~ \bar{x}}}k_{i}x^{v_{\cdot i}}\frac{e^{\alpha_{i}}}{2}[(v'_{\cdot i}-v_{\cdot i})^{\top}\triangledown f(x)]^2}_{H_{1}(x)}\\
	    &\underbrace{-\varlimsup\limits_{x\rightarrow\bar{x}}\sum_{\substack{\mathrm{supp}~ v_{\cdot i}\nsubseteq \mathrm{supp} ~\bar{x}\\\mathrm{supp}~v'_{\cdot i}\subseteq \mathrm{supp}~ \bar{x}}}k_{i}x^{v_{\cdot i}}\frac{e^{\alpha_{i}}}{2}[(v'_{\cdot i}-v_{\cdot i})^{\top}\triangledown f(x)]^2}_{H_{2}(x)},
	     \end{align*}
i.e., $$F_1(x)+F_3(x)\leq H_1(x)+H_{2}(x).$$ Obviously, $H_1(x)\leq 0$ and $H_{2}(x)\leq 0$.	     
 
We are only concerned with $H_1(x)$ to follow the proof. Denote $W=\text{supp}^c~ \bar{x}$, which is a semilocking set. Naturally, under the conditions of $\mathrm{supp}~v'_{\cdot i}\subseteq \mathrm{supp} ~\bar{x}$ and $ \mathrm{supp}~v_{\cdot i} \subseteq \mathrm{supp}~\bar{x}$, the $j$th component of $v'_{\cdot i}-v_{\cdot i}$ is zero for all $j\in W$. For simplicity but without loss of generality, we assume $W=\{m+1,~\cdots,~n\}$ where $m<n$. Then under the mentioned conditions, we have 
\begin{equation*}
	    (v'_{\cdot i}-v_{\cdot i})^{\top}\triangledown f(x)=\sum_{j=1}^{m}(v'_{ji}-v_{ji})\triangledown_{j}f(x).
	    \end{equation*}
Since $\Delta f(x)$ is a diagonal matrix, i.e.,
	    \begin{equation*}
	    \frac{\partial^2 f}{\partial x_{j_{1}}\partial x_{j_{2}}}=0=\frac{\partial \triangledown_{j_1} f}{\partial x_{j_{2}}},
	    \end{equation*}   
$\triangledown_{j_{1}}f(x)$ does not contain $x_{j_{2}}$ if $j_1\in W$ while $j_2\notin W$. By letting $x_\bot$ and $\bar{x}_\bot$ denote the first $m$ components of $x$ and $\bar{x}$, respectively, we thus get
	   \begin{align*}
	   H_{1}(x)
	   &=-\varlimsup\limits_{x_\bot\rightarrow \bar{x}_\bot}\sum_{\substack{\mathrm{supp} ~v_{\cdot i}\subseteq \mathrm{supp}~ \bar{x}\\\mathrm{supp}~v'_{\cdot i}\subseteq \mathrm{supp}~ \bar{x}} }k_{i}x_\bot^{v_{\cdot i}}\frac{e^{\alpha_{i}}}{2}\left(\sum_{j=1}^{m}(v'_{ji}-v_{ji}) \triangledown_{j} f(x_\bot)\right)^2\\
	   &=-\sum_{\substack{\mathrm{supp} ~v_{\cdot i}\subseteq \mathrm{supp}~ \bar{x}\\\mathrm{supp}~v'_{\cdot i}\subseteq \mathrm{supp}~ \bar{x}}}k_{i}\bar{x}_\bot^{v_{\cdot i}}\frac{e^{\alpha_{i}}}{2}\left(\sum_{j=1}^{m}(v'_{ji}-v_{ji})\triangledown_{j} f(\bar{x}_\bot)\right)^2.
	   \end{align*}
	   
If $H_1(x)=0$, which implies $\triangledown f(\bar{x}_\bot)\perp (v'_{\cdot i}-v_{\cdot i})$ for all $\mathrm{supp}~v_{\cdot i}\subseteq\mathrm{supp}~\bar{x}$ and $\mathrm{supp}~v'_{\cdot i}\subseteq\mathrm{supp}~\bar{x}$, we obtain 	   
 \begin{equation}\label{A.5}
	     \lim_{x\rightarrow \bar{x}}\sum_{\substack{\mathrm{supp} ~v_{\cdot i}\subseteq \mathrm{supp}~ \bar{x}\\\mathrm{supp}~v'_{\cdot i}\subseteq \mathrm{supp}~ \bar{x}}}k_{i}x^{v_{\cdot i}}=\lim_{x\rightarrow\bar{x}}\sum_{\substack{\mathrm{supp} ~v_{\cdot i}\subseteq \mathrm{supp}~ \bar{x}\\\mathrm{supp}~v'_{\cdot i}\subseteq \mathrm{supp}~ \bar{x}}}k_{i}x^{v_{\cdot i}}\exp[(v'_{\cdot i}-v_{\cdot i})^{\top}\triangledown f(x)]
	   \end{equation}	    
Combining (A.2), (A.4) and the boundary condition (a) (\ref{eq:pdea}), we further have 	    
\begin{equation}\label{eq:4.12}
    \lim_{x\rightarrow \bar{x}}\sum_{\substack{\mathrm{supp} ~v_{\cdot i}\nsubseteq \mathrm{supp}~ \bar{x}\\\mathrm{supp}~v'_{\cdot i}\subseteq \mathrm{supp}~ \bar{x}}}k_{i}x^{v_{\cdot i}}\exp [(v'_{\cdot i}-v_{\cdot i})^{\top}\triangledown f(x)]=0.
\end{equation}
This result also applies to other semilocking boundary non-equilibrium points, denoted by $\bar{x}^*$, satisfying $\mathrm{supp}~ \bar{x}=\mathrm{supp}~ \bar{x}^*$. The main reason is that the change of $\exp\{(v'_{\cdot i}-v_{\cdot i})^{\top}\triangledown f(x)\}$ from $x\rightarrow \bar{x}$ to $x\rightarrow \bar{x}^*$ in (\ref{eq:4.12}) will not change the limit value, which can be easily verified when we make the following decomposition  
$$(v'_{.i}-v_{.i})^{\top}\triangledown f(x)=\sum_{j=1}^{m} (v'_{ji}-v_{ji})\triangledown_{j} f(x_{j})+\sum_{j=m+1}^{n}(v'_{ji}-v_{ji})\triangledown_{j}f(x_{j}).$$
Clearly, the second term in the right hand of the above equation will keep unchanged from $x\rightarrow \bar{x}$ to $x\rightarrow \bar{x}^*$ since the last $n-m$ entries of $\bar{x}$ and $\bar{x}^*$ are the same. Although the first term will change, both $\bar{x_{j}}^{*}$ and $\bar{x_{j}}$ are specific positive bounded values for $j=1,\cdots,m$, which will lead to $f(\bar{x_{j}})$ and $f(\bar{x_{j}}^{*})$ bounded. Therefore, the first term will keep bounded for any boundary points satisfying $\mathrm{supp}~ \bar{x}=\mathrm{supp}~ \bar{x}^*$. Namely, (\ref{eq:4.12}) holds for any $\bar{x}^{*}$ with $\mathrm{supp}^{c}~\bar{x}^{*}=W$. We thus can generalize (\ref{A.5}) into
\begin{equation}\label{eq:4.14}
    \lim_{\substack{x_{j}>0,j=1,\cdots, m\\x_{j}\rightarrow 0, j=m+1,\cdots,n}}\sum_{\substack{\mathrm{supp} ~v_{\cdot i}\subseteq \mathrm{supp}~ \bar{x}\\\mathrm{supp}~v'_{\cdot i}\subseteq \mathrm{supp}~ \bar{x}}}k_{i}x^{v_{\cdot i}}=\lim_{\substack{x_{j}>0,j=1,\cdots, m\\x_{j}\rightarrow 0, j=m+1,\cdots,n}}\sum_{\substack{\mathrm{supp} ~v_{\cdot i}\subseteq \mathrm{supp}~ \bar{x}\\\mathrm{supp}~v'_{\cdot i}\subseteq \mathrm{supp}~ \bar{x}}}k_{i}x^{v_{\cdot i}}\exp[(v'_{\cdot i}-v_{\cdot i})^{\top}\triangledown f(x)],
\end{equation}
which is actually equivalent to 
\begin{equation}\label{eq:4.15}
   \sum_{\substack{\mathrm{supp} ~v_{\cdot i}\subseteq \mathrm{supp}~ \bar{x}\\\mathrm{supp}~v'_{\cdot i}\subseteq \mathrm{supp}~ \bar{x}}}k_{i}x_\bot^{\tilde{v}_{\cdot i}}=\sum_{\substack{\mathrm{supp} ~v_{\cdot i}\subseteq \mathrm{supp}~ \bar{x}\\\mathrm{supp}~v'_{\cdot i}\subseteq \mathrm{supp}~ \bar{x}}}k_{i}x_\bot^{\tilde{v}_{\cdot i}}\exp[(\tilde{v}'_{\cdot i}-\tilde{v}_{\cdot i})^{\top}\triangledown \tilde{f}(x)].
\end{equation}
Here, $\tilde{v}_{\cdot i}$, $\tilde{v}'_{\cdot i}$ and $\triangledown \tilde{f}(x)$ are the first $m$ entries of $v_{\cdot i},~v'_{\cdot i}$ and $\triangledown f(x)$, respectively. Now we define another MAS $(\mathcal{\tilde{S},\tilde{C},\tilde{R}},\tilde{k})$ with $\mathcal{\tilde{S}}=\{S_{1},\cdots,S_{m}\},~\mathcal{\tilde{R}}=\{v_{\cdot i}\to v'_{\cdot i}|~\mathrm{supp}~v_{\cdot i},\mathrm{supp}~v'_{\cdot i}\subseteq \{1,\cdots,m\}$\}, $\mathcal{\tilde{C}}=\{\text{complexes~in}~\mathcal{\tilde{R}}\}$ and $\tilde{k}$ as the projection of $k\in \mathbb{R}^{\mathcal{R}}$ into $\mathbb{R}^{\mathcal{\tilde{R}}}$. Then, according to the new MAS, (\ref{eq:4.14}) may be rewritten as
\begin{equation}
    \sum_{\mathcal{\tilde{R}}}\tilde{k}_{i}x_\bot^{\tilde{v}_{\cdot i}}=\sum_{\mathcal{\tilde{R}}}\tilde{k}_{i}x_\bot^{\tilde{v}_{\cdot i}}\exp\{(\tilde{v}'_{\cdot i}-\tilde{v}_{\cdot i})^{\top}\triangledown \tilde{f}(x)\}.
\end{equation}
Note that $\triangle f(x)$ is a diagonal matrix, so $\triangledown \tilde{f}(x)$ is independent of $x_j$ when $j=m+1,\cdots,n$. Denote the function whose gradient with respect to $x_\bot$ is $\triangledown \tilde{f}(x)$ by $\tilde{f}(x)$, then $\tilde{f}(x)$ is only related to  $x_\bot$. Therefore, $\tilde{f}(x)$ is a solution of (\ref{eq:4.14}), i.e., the Lyapunov Function PDE of the MAS $(\mathcal{\tilde{S},\tilde{C},\tilde{R}},\tilde{k})$. Also, note that $\triangle \tilde{f}(x)$ is positive definite, which together with the assumption of $H_1(x)=0$ suggests that $\bar{x}_\bot$ is a positive equilibrium point in $(\mathcal{\tilde{S},\tilde{C},\tilde{R}},\tilde{k})$ (based on $\textbf{Proposition \ref{pro:dissipation}}$). Therefore, $\bar{x}$ is a boundary equilibrium point in the MAS $(\mathcal{S,C,R},k)$, which is contradict to the condition that $\bar{x}$ is a boundary non-equilibrium point. We thus have $H_1(x)<0$, i.e., $\exists~M>0$ such that $F_1(x)+F_3(x)<-M$. 

From (\ref{eq:3.3}) in the proof of \textbf{Lemma \ref{lm:3.2}}, we have $F_4(x)\leq 0$, so $\varlimsup\limits_{x\rightarrow \bar{x}}\dot{f}(x)\leq -M$. Further from \textbf{Lemma \ref{lm:4.1}}, we get that $\bar{x}$ is not an $\omega$-$limit$ point. $\Box$
	   \end{appendix}

\end{document}